\documentclass[12pt]{article}

\usepackage{amsmath,a4wide,amssymb,amsfonts,amsbsy,xy,latexsym,epsf,pstricks,times}
\xyoption{all}

\usepackage[latin1]{inputenc}
\usepackage[french]{babel}

\let\set\mathbb

\newenvironment{preuve}{\medbreak \noindent {\bf Preuve~---}}
                       {\hfill $\square$ \medbreak}

\def\QQ{{\set Q}}

\def\Spec{\mathop{\rm{Spec}}\nolimits }
\def\Id{{\rm Id}}
\def\Qb{{\overline{\QQ}}}
\def\bQ{{\overline{\QQ}}}

\def\GK{{\Gamma_K}}

\def\pgcd{\mathop{\rm{pgcd}}\nolimits }

\def\Spec{\mathop{\rm{Spec}}\nolimits }
\def\Proj{\mathop{\rm{Proj}}\nolimits }

\def\Gal{\mathop{\rm{Gal}}\nolimits }
\def\Aut{\mathop{\rm{Aut}}\nolimits }

\def\cI{{\cal I}}

\def\cV{{\cal V}}

\def\cP{{\cal P}}

\def\cE{{\cal E}}
\def\cQ{{\cal Q}}
\def\cJ{{\cal J}}
\def\cK{{\mathcal K}}

\def\cM{{\cal M}}

\def\cO{{\cal O}}

\def\agot{{\mathfrak a}}
\def\bagot{{\bar {\mathfrak a}}}

\def\bgot{{\mathfrak b}}

\newtheorem{theoreme}{Th{\'e}or{\`e}me}[section]
\newtheorem{lem}[theoreme]{Lemme}
\newtheorem{definition}[theoreme]{D{\'e}finition}

\newtheorem{corollaire}[theoreme]{Corollaire}

\newtheorem{remarque}{Remarque}

\def\cC{{\cal C}}
\def\cK{{\cal K}}
\def\cJ{{\cal J}}
\def\cU{{\cal U}}
\def\cG{{\cal G}}
\def\PP{{\set P}}
\def\PU{{{\set P}^1}}

\def\cM{{\cal M}}

\def\cS{{\cal S}}

\def\cB{{\cal B}}

\def\cD{{\cal D}}
\def\cF{{\cal F}}

\def\NN{{\set N}}
\def\cO{{\cal O}}

\def\ZZ{{\set Z}}

\def\KK{{ K\,}}

\def\LL{{ L\,}}

\def\CC{{\set C}}
\def\AA{{\set A}}
\def\PP{{\set P}}

\begin{document}
\author{Jean-Marc Couveignes et Emmanuel Hallouin\thanks{Institut de Math{\'e}matiques de Toulouse,
    Universit{\'e} de Toulouse et CNRS, Universit{\'e}  de Toulouse 2
le Mirail, 5 all{\'e}es
    Antonio Machado 31058 Toulouse c{\'e}dex 9}}
\title{Sur le corps des modules  de certaines vari{\'e}t{\'e}s}

\maketitle

\begin{abstract}
Nous associons {\`a} tout rev{\^e}tement de courbes
des vari{\'e}t{\'e}s ayant m{\^e}me corps des modules et m{\^e}mes corps
de d{\'e}finition que ce rev{\^e}tement. Nous en d{\'e}duisons  des exemples
de courbes qui ont $\QQ$ pour corps des modules, admettent
des mod{\`e}les sur tous les compl{\'e}t{\'e}s de $\QQ$ mais
pas de mod{\`e}le sur $\QQ$.
\end{abstract}

\tableofcontents

\section{Introduction}\label{section:introduction}

Le propos de cet article est de construire des obstructions {\`a} la
descente dans la cat{\'e}gorie des vari{\'e}t{\'e}s, c'est-{\`a}-dire
d'exhiber,
 par
exemple, des
vari{\'e}t{\'e}s de corps des modules $\QQ$ n'ayant pas de mod{\`e}le sur~$\QQ$. 

Partant d'une obstruction {\`a} la descente ({\'e}ventuellement  globale comme
dans~\cite{roscouveignes}) dans la
cat{\'e}gorie des rev{\^e}tements de courbes, on souhaite produire des
obstructions
{\`a} la descente dans d'autres cat{\'e}gories. La cat{\'e}gorie ultime sera celle des
courbes lisses. On proc{\'e}dera par {\'e}tapes : passant par la cat{\'e}gorie des
surfaces quasi-projectives, puis projectives.

Peu d'exemples d'obstructions sont connus pour les vari{\'e}t{\'e}s. Mestre a donn{\'e} des exemples
d'obstructions locales pour des courbes hyperelliptiques dans
\cite{Mestre}.

On montre ici le

\begin{theoreme}\label{theoreme:courbex}
Il existe une $\overline{\QQ}$-courbe lisse projective, de corps de module~$\QQ$ et d{\'e}finie sur tous les
compl{\'e}t{\'e}s de~$\QQ$ mais pas sur~$\QQ$ lui m{\^e}me.
\end{theoreme}

\subsection{Quelques cat{\'e}gories  munies d'une action de Galois} \label{subsection:5cat}

Soit~$K$ un corps de caract{\'e}ristique nulle.
On note~$K^s$  sa cl{\^o}ture s{\'e}parable et~$\Gamma_K$ son groupe de
Galois absolu. Une  $K$-vari{\'e}t{\'e} est par d{\'e}faut suppos{\'e}e
r{\'e}duite, lisse  et g{\'e}om{\'e}triquement irr{\'e}duc\-ti\-ble. 
Nous serons amen{\'e}s {\`a} consid{\'e}rer diff{\'e}rentes cat{\'e}gories munies d'une action
fonctorielle
de $\GK$.
Il s'agit premi{\`e}rement de la cat{\'e}gorie des vari{\'e}t{\'e}s
alg{\'e}briques, r{\'e}duites, lisses, et irr{\'e}ductibles
sur $K^s$. On consid{\`e}re aussi les sous cat{\'e}gories pleines form{\'e}es  des
courbes lisses projectives sur $K^s$, des surfaces quasi-projectives lisses
sur $K^s$, des
surfaces
projectives normales sur $K^s$. Soit $\cB$ une $K$-vari{\'e}t{\'e} lisse
r{\'e}duite et g{\'e}om{\'e}triquement irr{\'e}ductible.
On s'int{\'e}ressera {\`a} la cat{\'e}gorie des $K^s$-rev{\^e}tements de 
$\cB$ et {\`a} la sous-cat{\'e}gorie pleine form{\'e}e des rev{\^e}tements {\'e}tales.

Si $\cO \rightarrow \Spec (K^s)$
est un objet d'une des  cat{\'e}gories ci-dessus
et si $\sigma \in \GK$, l'objet   ${}^\sigma\!\cO\rightarrow
\Spec (K^s)$
est d{\'e}fini par composition  $\cO\rightarrow  \Spec(K^s)
\stackrel{\scriptscriptstyle \Spec (\sigma ^{-1} ) }{\longrightarrow} \Spec
(K^s)$. Cela d{\'e}finit un foncteur covariant, encore not{\'e} $\sigma$.
 Le stabilisateur dans $\GK$ de la
classe de $K^s$-isomorphisme d'un objet est un
sous-groupe
d'indice fini. Le sous-corps de $K^s$ fix{\'e} par ce sous-groupe est
appel{\'e} corps
des modules de l'objet. 
Si $L$ est
un corps tel que $K\subset L\subset K^s$, on dit  que $\cO$ est d{\'e}fini sur $L$
 s'il   est $K^s$-isomorphe au tir{\'e} en arri{\`e}re d'un $L$-objet $\cO'\rightarrow \Spec (L)$
 par $\Spec (K^s)\rightarrow \Spec (L)$. 
On dit aussi  que $L$ est un  corps de d{\'e}finition
de $\cO$. Le
corps
des modules est contenu dans  tous  les corps de d{\'e}finition. 

\subsection{Descente de Weil}

Soit $K$ un corps de caract{\'e}ristique nulle
 et $K^s$ une cl{\^o}ture alg{\'e}brique.
Soit $L\subset K^s$ une extension finie de $K$.
Soit $\Theta$ l'ensemble des $K$-plongements de $L$ dans $K^s$.
C'est un ensemble fini de cardinal le degr{\'e} de $L$ sur $K$. Il est
muni d'une action {\`a} gauche du groupe  de Galois absolu $\GK$ de
$K$.
Pour $\omega \in \GK$ et $\sigma \in \Theta$ cette action est
not{\'e}e
sans fa{\c c}ons $\omega\sigma\in \Theta$.  \`A l'inclusion
$L\subset K^s$  correspond un
{\'e}l{\'e}ment particulier de $\Theta$ que l'on note $\Id$. 
On a une surjection $\GK\rightarrow \Theta$ d{\'e}finie par
$\omega \rightarrow \omega\Id$ et  on notera abusivement $\omega \Id
=\omega \in \Theta$.

Soit $\cQ\rightarrow \Spec (K^s)$ un objet de l'une des cat{\'e}gories
pr{\'e}c{\'e}dentes. 
On suppose que $K$ est le corps des modules de $\cQ$. 
On suppose que $\cQ$ est d{\'e}fini sur $L$.
Il existe donc un objet $\cO\rightarrow \Spec (L)$ tel que $\cQ$ soit
$K^s$-isomorphe   {\`a} 
${}^{\Id}\cO$.

On appelle $\langle \cO, L/K \rangle $ la sous-cat{\'e}gorie  pleine
dont les objets
sont les conjugu{\'e}s ${}^\sigma  \cO$   de $\cO$ avec $\sigma $
parcourant $\Theta$. Ici, l'objet ${}^\sigma  \cO$ est le tir{\'e} en
arri{\`e}re
de $\cO\rightarrow \Spec (L)$ le long de 
$\Spec(\sigma) : \Spec(K^s)\rightarrow \Spec(L)$.


Notons que $\langle \cO, L / K  \rangle $ 
est un groupo{\"\i}de fini dont le nombre d'objets
est le degr{\'e} de $L$ sur $K$. Ce groupo{\"\i}de est connexe puisque $K$
 est  le corps
des modules de $\cO$.
Ce groupo{\"\i}de est muni d'une action fonctorielle du groupe
de Galois absolu de $K$. Si $\cO/L$ est d{\'e}fini sur $L$ et si $M\subset K^s$
est une extension finie de $L$ de degr{\'e} $d$, alors $\langle\cO\otimes_LM,M/K\rangle$ est form{\'e} de
$d$  copies de $\langle \cO,L/K\rangle $ connect{\'e}es avec les identification triviales.

Weil montre que le groupo{\"\i}de galoisien abstrait $<\cO,L/K>$
contient toutes
les informations utiles {\`a} la descente de $L$ {\`a} $K$.

\begin{theoreme}[A. Weil]
L'objet $\cO$  est $L$-isomorphe {\`a} un $K$-objet
de sa cat{\'e}gorie si et seulement s'il existe une famille $(I_{\sigma, \tau})_{(\sigma, \tau)\in \Theta}$  de
morphismes dans la cat{\'e}gorie   $\langle\cO, L/K \rangle $ satisfaisant les conditions
suivantes. 

\begin{itemize}
\item Cette famille est   index{\'e}e par les
couples de plongements $(\sigma, \tau)\in \Theta$ et $I_{\sigma, \tau}
:  {}^\tau  \cO\rightarrow {}^\sigma  \cO$. 
\item $I_{\sigma, \tau}I_{\tau, \mu}=I_{\sigma, \mu}$ pour tous
  $\sigma, \tau, \mu \in \Theta$.
\item Pour tout $\omega\in \GK$ on a ${}^\omega I_{\sigma, \tau}=I_{\omega\sigma, \omega\tau}$.
\end{itemize}

\end{theoreme}

La famille de morphismes $(I_{\sigma, \tau})_{(\sigma, \tau)\in
  \Theta}$  est appel{\'e}e {\it donn{\'e}e de descente}
 de $L$ {\`a} $K$ pour $\cO$.

En somme, Weil montre que, pour les cat{\'e}gories concern{\'e}es,
les seules obstructions {\`a} la descente sont combinatoires et
qu'elles sont exprim{\'e}es en termes du groupo{\"\i}de galoisien abstrait
$\langle \cO, L/K \rangle $.
En particulier, soient  deux $L$-objets $\cO$ et $\cP$ 
de  corps des modules  $K$. Soient 
$\langle \cO, L/K\rangle $ et $\langle\cP, L/K\rangle $  
les
groupo{\"\i}des pour la descente {\`a} $K$. S'il existe un
isomorphisme de groupo{\"\i}des $\GK$-{\'e}quivariant entre $\langle \cO, L/K
\rangle $ et $\langle\cP, L/K \rangle $, alors $K$
est corps de d{\'e}finition de $\cO$ si et seulement s'il est corps de
d{\'e}finition
de $\cP$. Noter que $\cO$ et $\cP$ ne sont pas, ici, n{\'e}cessairement objets
de la m{\^e}me cat{\'e}gorie : seuls importent les deux groupo{\"\i}des 
galoisiens associ{\'e}s.

Mais on a mieux : la seule existence d'un morphisme entre ces deux
groupo{\"\i}des fournit une information.

\begin{theoreme}\label{theoreme:criterech}
Soit $K$ un corps de caract{\'e}ristique nulle.
Soit $\cO$ un objet de l'une des  cat{\'e}gories ci-dessus.
 Soit
$\cP$ un objet de l'une des  cat{\'e}gories ci-dessus (pas
n{\'e}cessairement celle de $\cO$).

Soit  $L\subset K^s$ une extension finie de $K$. On suppose que
$\cO$ et $\cP$ sont d{\'e}finis sur $L$ et que $K$ est
le corps des modules de  $\cO$ et $\cP$.  Soient $\langle\cO, L/K
\rangle $  et $\langle\cP, L/K \rangle $ les
groupo{\"\i}des associ{\'e}s {\`a} $\cO$ et $\cP$ pour la descente de $L$ {\`a}
$K$.

Supposons qu'il existe un morphisme (foncteur covariant)  de groupo{\"\i}de
$\GK$-{\'e}quivariant 
de $\langle\cO, L/K \rangle $ vers
$\langle\cP, L/K \rangle $ qui envoie l'objet $\cO$ sur l'objet $\cP$. 

Alors 
si $K$ est corps de d{\'e}finition de $\cO$, il est aussi corps de
d{\'e}finition de $\cP$.
\end{theoreme}

C'est {\'e}vident : il suffit d'appliquer le foncteur en question  {\`a}
une donn{\'e}e de descente $(I_{\sigma, \tau})_{(\sigma, \tau)\in \Theta}$
de $L$ {\`a} $K$ pour $\cO$. On obtient une donn{\'e}e de descente de
$L$ {\`a} $K$ pour $\cP$ car le foncteur  est un morphisme
$\GK$-{\'e}quivariant
de groupo{\"\i}des.

\begin{remarque}
Le groupo\"ide muni d'une action de Galois qui appara{\^\i}t dans les th{\'e}or{\`e}mes
pr{\'e}c{\'e}\-dents est tr{\`e}s proche de la gerbe des mod{\`e}les introduite  par
Giraud dans \cite{giraud} et utilis{\'e}e par D{\`e}bes et Douai dans \cite{DD}
dans le cas des rev{\^e}tements. La formulation de Weil fait appara{\^\i}tre 
seulement la famille finie des conjugu{\'e}s d'un mod{\`e}le donn{\'e}
(soit des sections locales de la gerbe sur un recouvrement fini de
$\Spec(K)$). 
Les conditions de Weil 
sont {\'e}quivalentes {\`a} l'existence 
d'un groupo\"ide galoisien $\cG$ dont les objets sont les objets
de $\langle  \cO, \LL/\KK  \rangle$ plus un objet $O$, et qui satisfait les
propri{\'e}t{\'e}s suivantes :
\begin{itemize}
\item Il existe un morphisme $I$ dans $\cG$ entre $O$ et $\cO$.
\item Le groupo\"ide $\langle  \cO, \LL/\KK \rangle$  est une sous-cat{\'e}gorie pleine
de  $\cG$.
\item L'action du groupe de Galois absolu de $\KK$ sur $\langle  \cO, \LL/ \KK  \rangle$
s'{\'e}tend en une action fonctorielle sur $\cG$ qui fixe $O$.
\item Le morphisme  $I$ est fix{\'e} par le groupe de Galois absolu
de $\LL$ vu comme  sous-groupe du groupe de Galois absolu
de $\KK$.
\end{itemize}

\end{remarque}
\subsection{Strat{\'e}gie}

Les obstructions {\`a} la descente construites dans ce travail 
d{\'e}rivent 
de celles construites par   Ros
et  Couveignes  dans la cat{\'e}gorie des rev{\^e}tements
de courbes :

\begin{theoreme}[cf.\cite{roscouveignes}, Corollaire~2] \label{theoreme:catrevetement}
Il existe un $\overline{\QQ}$-rev{\^e}tement ramifi{\'e} de~$\PP^1_\QQ$ dont le corps des modules
est~$\QQ$, qui admet des mod{\`e}les sur tous les compl{\'e}t{\'e}s de~$\QQ$ mais pas de mod{\`e}le sur~$\QQ$.
\end{theoreme}

Afin d'en d{\'e}duire des exemples d'obstructions {\`a} la descente dans
d'autres cat{\'e}gories, nous pr{\'e}sentons plusieurs proc{\'e}d{\'e}s
permettant de passer de la cat{\'e}gorie des rev{\^e}tements de courbes aux 
autres cat{\'e}gories pr{\'e}sent{\'e}es
dans la section~\ref{subsection:5cat}, tout en conservant les propri{\'e}t{\'e}s de modules et de
d{\'e}finition. Pour ce faire, nous
utiliserons plusieurs fois le th{\'e}or{\`e}me~\ref{theoreme:criterech} pour passer d'une
cat{\'e}gorie {\`a} une autre.

Dans la section \ref{section:sansauto}, on part d'un
$K^s$-rev{\^e}tement $\varphi : \cC\rightarrow \cB\otimes_K K^s$  d'une $K$-courbe $\cB$ et on construit une
$K$-courbe 
$\cB'$ sans $K^s$-automorphisme non-trivial et un $K^s$-rev{\^e}tement 
$\cC'\rightarrow \cB'\otimes_KK^s$ ayant m{\^e}me corps des modules et m{\^e}mes corps de
d{\'e}finition que $\varphi$.

Dans la section \ref{section:surfacesaffines} on part d'un
$K^s$-rev{\^e}tement $\varphi : \cC\rightarrow \cB\otimes_KK^s$  d'une
$K$-courbe $\cB$ et on construit
 une $K^s$-surface quasi-projective  ayant m{\^e}me corps des modules et m{\^e}mes
corps de d{\'e}finition que ce rev{\^e}tement. Lorsque la base $\cB$
de $\varphi$ 
n'a pas de $K^s$-automorphisme non-trivial, la surface en question est
l'ouvert compl{\'e}mentaire du graphe du rev{\^e}tement. Le r{\'e}sultat
de la section \ref{section:sansauto} permet de se ramener {\`a} ce cas.


Dans la section \ref{section:surfacesnormalesprojectives} nous partons d'un  $K^s$-rev{\^e}tement 
d'une $K$-courbe
alg{\'e}brique sans $K^s$-automorphisme. Nous supposons que le corps des modules de ce rev{\^e}tement est
$K$ et nous lui associons une $K^s$-surface projective normale, de corps des
modules $K$,    ayant les  m{\^e}mes
corps de d{\'e}finition que ce rev{\^e}tement.  Cette surface projective
est construite comme  rev{\^e}tement
fortement ramifi{\'e} le long du graphe du rev{\^e}tement.

Enfin, dans la section
\ref{section:courbes}  nous construisons une $K^s$-courbe projective, de corps
des modules $K$, 
ayant les  m{\^e}mes
corps de d{\'e}finition que le rev{\^e}tement initial.
 C'est une courbe trac{\'e}e sur la surface pr{\'e}c{\'e}dente, et obtenue
par d{\'e}formation d'une courbe stable choisie pour que son groupe
d'automorphismes soit le m{\^e}me que celui de ladite surface.

Tout au long de ce travail, notre principale pr{\'e}occupation sera
de contr{\^o}ler  le groupe des automorphismes de chacun des objets
que nous construirons.

\section{Annulation du groupe d'automorphismes de la base}\label{section:sansauto}
Dans cette section, on part d'un corps $K$ de caract{\'e}ristique nulle,
d'une extension alg{\'e}brique
 $L\subset K^s$ de $K$, d'une $K$-courbe~$\mathcal{B}$
projective, lisse, g{\'e}om{\'e}triquement irr{\'e}ductible,  et d'un
$L$-rev{\^e}tement~$\varphi : \mathcal{C} \to \mathcal{B} \otimes_K L$ ayant~$K$ pour corps des modules.
On veut produire d'autres  rev{\^e}tements, ayant m{\^e}me  corps des modules
 et m{\^e}mes  corps de d{\'e}finition, et jouissant de propri{\'e}t{\'e}s
 suppl{\'e}mentaires.
 En particulier, on souhaite montrer que l'on peut
supposer que la base ne poss{\`e}de pas de $K^s$-automorphisme non trivial.

On commence par montrer que le degr{\'e} du rev{\^e}tement peut {\^e}tre multipli{\'e} par n'importe quel premier
{\'e}tranger au degr{\'e} de d{\'e}part.

\begin{lem} \label{lem:granddegre}
Soit~$\mathcal{B}$ une $K$-courbe projective lisse g{\'e}om{\'e}triquement irr{\'e}ductible
et~$\varphi : \mathcal{C} \to \mathcal{B} \otimes_K L$ un $L$-rev{\^e}tement
de degr{\'e}~$d$. Pour tout
premier~$p$ {\'e}tranger  {\`a}~$d$, il existe~$\varphi' : \mathcal{C}' \to \mathcal{B} \otimes_K L$
un $L$-rev{\^e}tement de degr{\'e}~$pd$ ayant m{\^e}me corps des modules et m{\^e}mes corps de d{\'e}finition que~$\varphi$.
\end{lem}

\begin{preuve}
Soit~$f \in K(\mathcal{B})$ une fonction non constante, dont le diviseur est simple et disjoint
du lieu de ramification de~$\varphi$. L'{\'e}quation~$Y^p = f$ d{\'e}finit une extension de degr{\'e}~$p$
de~$K(\mathcal{B})$ dont on note~$\mathcal{B}'$ le mod{\`e}le projectif lisse. C'est une $K$-courbe,
projective, lisse, g{\'e}om{\'e}triquement irr{\'e}ductible qui rev{\^e}t~$\mathcal{B}$ par un rev{\^e}tement~$\nu$
de degr{\'e}~$p$ (et galoisien sur~$K^s$). Les extensions~$L(\mathcal{C})$
et~$L(\mathcal{B}')$ sont lin{\'e}airement disjointes sur~$L(\mathcal{B})$.
Soit~$\mathcal{C}'$ le mod{\`e}le projectif lisse associ{\'e} au compositum de~$L(\mathcal{C})$
et~$L(\mathcal{B}')$. C'est, par construction, une $L$-courbe qui rev{\^e}t~$\mathcal{B} \otimes_K L$
par un rev{\^e}tement~$\psi$ de degr{\'e}~$pd$~:
$$
\xymatrix@=15pt{
& \mathcal{C}' \ar@{->}[dl] \ar@{->}[dr] \ar@{->}[dd]^{\psi} \\
\mathcal{B}' \otimes_K L \ar@{->}[dr]_{\nu} & & \mathcal{C} \ar@{->}[dl]^{\varphi} \\
& \mathcal{B} \otimes_K L
}
$$
De nouveau par construction, il est facile de voir que tout corps de d{\'e}finition de~$\varphi$ est
aussi un corps de d{\'e}finition pour~$\psi$, et que  le corps des modules
de $\psi$ est inclus dans celui de $\varphi$. 

R{\'e}ciproquement, pour tout $\sigma \in \GK$, le
sous-rev{\^e}tement~${}^\sigma \varphi : {}^\sigma \mathcal{C} \to \mathcal{B} \otimes_K L$
du rev{\^e}tement~${}^\sigma \psi : {}^\sigma \mathcal{C}' \to \mathcal{B} \otimes_K L$ peut {\^e}tre
caract{\'e}ris{\'e} comme le sous-rev{\^e}tement maximal non ramifi{\'e} en
le support  de~$f$.
Il en r{\'e}sulte que pour tous~$\sigma,\tau \in \GK$ et tout $K^s$-isomorphisme~$I$
entre~${}^\tau \psi$ et~${}^\sigma \psi$, 
il existe un {\it unique} $K^s$-isomorphisme  $J : {}^\tau \cC \rightarrow {}^\sigma \cC$
tel que le diagramme

$$
\xymatrix@=15pt{
{}^\tau \cC' \ar@{->}[r]^I  \ar@{->}[d] & {}^\sigma \cC'  \ar@{->}[d] \\
{}^\tau \cC \ar@{->}[r]^J   & {}^\sigma \cC   \\
}
$$
\noindent commute.
 De plus $J$ est un $K^s$-isomorphisme entre ${}^\tau \varphi$ 
et ${}^\sigma \varphi$.
Cela permet
de montrer, d'une part, que le
corps des modules  de~$\varphi$ est contenu dans celui de~$\psi$. Compte tenu de ce qui pr{\'e}c{\`e}de ces
deux corps de modules sont {\'e}gaux. D'autre part, soit
$L' \subset K^s$, une extension
alg{\'e}brique de $M$. \`A  toute donn{\'e}e de descente  {\`a}~$L'$
pour~$\psi$, 
il correspond, par projection, une donn{\'e}e de
descente {\`a}~$L'$ pour~$\varphi$. Ainsi les rev{\^e}tements~$\varphi$ et~$\psi$ ont bien m{\^e}mes corps
de d{\'e}finition comme annonc{\'e}.
\end{preuve}

Ensuite, on montre que l'on peut choisir une courbe de base de genre~$\geq 2$.

\begin{lem} \label{lem:grandgenre}
Soit~$\mathcal{B}$ une $K$-courbe projective lisse g{\'e}om{\'e}triquement irr{\'e}ductible
et~$\varphi : \mathcal{C} \to \mathcal{B} \otimes_K L$ un $L$-rev{\^e}tement de degr{\'e}~$d$. Alors il
existe une $K$-courbe~$\mathcal{B'}$ lisse, projective, g{\'e}om{\'e}tri\-que\-ment irr{\'e}ductible, de genre au
moins~$2$ et
un $L$-rev{\^e}tement~$\varphi' : \mathcal{C}' \to \mathcal{B}' \otimes_K L$
ayant m{\^e}me corps de modules et m{\^e}mes corps de d{\'e}finition que celui de~$\varphi$.
\end{lem}

\begin{preuve}
On reprend la construction de la preuve pr{\'e}c{\'e}dente en supposant de plus que le degr{\'e} de~$f$
est au moins~$3$ et on s'int{\'e}resse au rev{\^e}tement~$\varphi' : \mathcal{C}' \to \mathcal{B}'\otimes_K L$.
L'hypoth{\`e}se sur le degr{\'e} de~$f$ permet, gr{\^a}ce {\`a} la formule de Hurwitz, de montrer que la
courbe~$\mathcal{B}'$ a un genre au moins {\'e}gal {\`a}~$2$. Avec des arguments proches de ceux
de la preuve pr{\'e}c{\'e}dente, on v{\'e}rifie que~$\varphi$ et~$\varphi'$ ont bien m{\^e}me corps des modules
et m{\^e}mes corps de d{\'e}finition.
\end{preuve}

Enfin, on prouve que la base peut ne pas avoir de $K^s$-automorphisme non trivial.

\begin{lem} \label{lem:sansauto}
Soit~$\mathcal{B}$ une $K$-courbe projective lisse g{\'e}om{\'e}triquement irr{\'e}ductible
et~$\varphi : \mathcal{C} \to \mathcal{B} \otimes_K L$ un $L$-rev{\^e}tement de degr{\'e}~$d$. Alors il
existe une $K$-courbe~$\mathcal{B'}$ lisse, projective,
g{\'e}om{\'e}trique\-ment irr{\'e}ductible, de genre au
moins~$2$ et telle que~$\mathcal{B'} \otimes_K K^s$ n'ait pas d'automorphisme non trivial
et il existe un $L$-rev{\^e}tement~$\varphi' : \mathcal{C}' \to \mathcal{B}' \otimes_K L$
ayant m{\^e}me corps de modules et m{\^e}mes corps de d{\'e}finition que~$\varphi$.
\end{lem}

\begin{preuve}
D'apr{\`e}s le lemme \ref{lem:grandgenre} 
 on peut supposer que le genre de $\cB$ est au moins $2$ et donc que 
 le groupe des $K^s$-automorphismes de $\cB$, not{\'e}~$A$, est fini.

Soit $p\ge 3$ un entier premier. Commen{\c c}ons par montrer qu'il existe une
fonction~$f \in K(\cB)$,
non-constante et non-singuli{\`e}re au-dessus de $2$, $-2$ et $\infty$,
de degr{\'e} au moins $2+4p(g(\cB)-1)+2p^2$,
telle que l'ensemble  $f^{-1}(\{-2,2\})$ ne soit invariant par aucun automorphisme non-trivial de $\cB\otimes_KK^s$,
et telle que l'ensemble des  valeurs singuli{\`e}res de $\varphi$ ne rencontre pas 
l'ensemble $f^{-1}(\{2,-2,\infty \})$.

En {\'e}vitant les noyaux des formes lin{\'e}aires~$\agot \pm {\rm id}$
pour~$\agot \in A \setminus \{{\rm id}\}$
(en nombre fini) du $K^s$-espace vectoriel~$K^s(\cB)$, on v{\'e}rifie qu'il existe~$f \in K^s(\cB)$ non
constante telle que~$\agot(f) \not= \pm f$ pour tout~$\agot \in A \setminus \{{\rm id}\}$.
On peut m{\^e}me supposer que~$f \in K(\cB)$.
Pour tout~$n \in \NN^*$ premier {\`a}~$\#A$ et tout~$\agot \in A \setminus \{{\rm id}\}$, on a
encore~$\agot(f^n) \not= \pm f^n$ (Sinon~$\agot(f^n) = \pm f^n$ donc~$\agot(f) = \zeta f$
avec~$\zeta^{2n} = 1$. Notant~$\alpha$ l'ordre de~$\agot$,
alors~$f = \agot^{\alpha}(f) = \zeta^{\alpha} f$
donc~$\zeta^{\alpha} = 1$. Ainsi~$\zeta = \pm 1$. Contradiction). Quitte {\`a} remplacer~$f$ par~$f^n$,
on peut supposer que~$\deg(f) \geq 2+4p(g(\cB)-1)+2p^2$.

Par construction, la fonction~$f^2$ n'a pas d'automorphisme non trivial
(i.e.~$\Aut_{K(f^2)}(K(\cB)) = \{{\rm id}\}$). D'apr{\`e}s le
lemme~\ref{lem:autofibres}, on en d{\'e}duit que presque toutes les fibres
de~$f^2$ sont non singuli{\`e}res et fix{\'e}es par aucun automorphisme non trivial de~$A$. En particulier,
il existe~$\lambda \in K^*$ telle que la
fibre de $f^2$ en $\lambda^2$ est non singuli{\`e}re, fix{\'e}e par aucun
automorphisme non-trivial 
de~$A$ et disjointe des valeurs singuli{\`e}res de~$\varphi$.
La fonction~$2f/\lambda$ satisfait toutes les contraintes souhait{\'e}es. On la note encore~$f$.

Cela {\'e}tant, l'{\'e}quation $X^p+X^{-p}-f=0$ d{\'e}finit une extension  de
$K(\cB)$, di{\'e}drale de groupe~$D_p$ apr{\`e}s extension des scalaires {\`a}  $K^s$.
Soit $\cB''$ le mod{\`e}le projectif lisse de ce corps de fonctions et 
$\mu  : \cB'' \rightarrow \cB$ le $K$-rev{\^e}tement de degr{\'e} $2p$ associ{\'e}.
Le groupe des $K^s$-automorphismes de $\mu$ est $D_p$.
L'ensemble des valeurs singuli{\`e}res de $\mu$ est $f^{-1}(\{2,-2,\infty \})$.
On note~$w$ l'automorphisme de $\cB''$ donn{\'e} par $w(X)=X^{-1}$
et  $\cB'$
le quotient de $\cB''$ par $w$.  C'est une $K$-courbe
projective et lisse et on a  un $K$-rev{\^e}tement
$\nu : \cB'\rightarrow \cB$ de degr{\'e} $p$. Ce rev{\^e}tement n'a pas d'automorphisme non trivial
car~$\langle w \rangle$ est auto-normalisateur dans~$D_p$.

Ramification disjointe oblige, les corps~$K^s(\cB'')$ et~$K^s(\cC)$ sont lin{\'e}airement disjoints
sur~$K^s(\cB)$. Les mod{\`e}les projectifs lisses des compositum de~$K^s(\cC)$ avec~$K^s(\cB')$
et~$K^s(\cB'')$, not{\'e}s~$\cC'$ et~$\cC''$ s'inscrivent dans le diagramme suivant~:
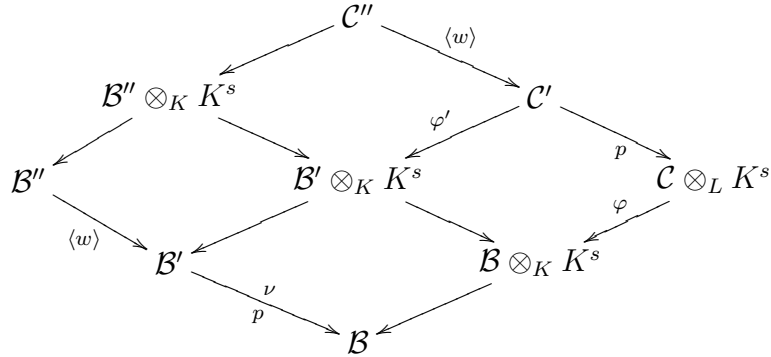
\begin{equation*}
\xymatrix@=15pt{
&&\cC''\ar@{->}[dl]\ar@{->}[dr]^{\langle w\rangle}& & \\
&\cB'' \otimes_KK^s  \ar@{->}[dl]  \ar@{->}[dr]&&\cC'  \ar@{->}[dl]_{\varphi'}\ar@{->}[dr]_p&\\
\cB'' \ar@{->}[dr]_{\langle w\rangle} &&\cB'\otimes_KK^s  \ar@{->}[dl]  \ar@{->}[dr]&&\cC\otimes_LK^s\ar@{->}[dl]_\varphi\\
&\cB' \ar@{->}[dr]^\nu_p &&\ar@{->}[dl]\cB\otimes_KK^s&  \\
&&\cB &}\label{diagramme:rev}
\end{equation*}
Le rev{\^e}tement~$\cC'' \to \cC$ est encore galoisien de groupe~$D_p$ et $\cC'' \to \cC'$
est encore de degr{\'e}~$2$.

Montrons que le rev{\^e}tement~$\varphi' : \cC' \to \cB'\otimes_KK^s$
satisfait les conditions du lemme.
D{\'e}j{\`a}, il a clairement m{\^e}me corps des modules et m{\^e}mes corps 
de d{\'e}finition que $\varphi$. 

Il suffit de montrer
 que $\cB'\otimes_KK^s$ n'a pas d'automorphisme non trivial. Soit $\bgot'$
un automorphisme de $\cB'\otimes_KK^s$. L'image $\cI$ de  $\nu\times (\nu\circ \bgot') : \cB' \otimes_KK^s
\rightarrow (\cB\times\cB)\otimes_KK^s$ est coinc{\'e}e entre $\cB'\otimes_KK^s$ et $\cB\otimes_KK^s$. Le
bidegr{\'e}  de $\cI$ est $\le (p,p)$ et donc son genre arithm{\'e}tique virtuel
est  au plus $1+2p(g(\cB)-1)+p^2$ d'apr{\`e}s le lemme \ref{lem:adjonction}. Si $\cI$
n'est pas isomorphe {\`a} $\cB\otimes_KK^s$
alors $\cI$ est birationnellement {\'e}quivalent {\`a} $\cB'\otimes_KK^s$ et son genre
g{\'e}om{\'e}trique est $>\frac{1}{4}\deg (f)p\ge 1+2p(g(\cB)-1)+p^2$ par
la formule
de Hurwitz. Contradiction.  Donc $\cI$ est une correspondance  de bidegr{\'e} $(1,1)$
qui d{\'e}finit un automorphisme $\bgot$ de $\cB\otimes_KK^s$
 tel que $\bgot\circ \nu =\nu \circ \bgot'$. Un tel automorphisme conserve 
les donn{\'e}es de ramification de $\nu$, de sa cl{\^o}ture galoisienne
 $\mu$ et aussi de l'unique sous-rev{\^e}tement de degr{\'e} deux de $\mu$.
Comme  ce dernier rev{\^e}tement est ramifi{\'e} en $f^{-1}(\{-2,2\})$ on en d{\'e}duit
que $\bgot =\Id$ puis que  $\bgot'$ est un $K^s$-automorphisme du rev{\^e}tement
$\nu$. Comme $\nu$ n'a pas d'automorphisme, on a bien~$\bgot' = {\rm id}$.
\end{preuve}

\section{Surfaces quasi-projectives}\label{section:surfacesaffines}


Soit~$K$ un corps de caract{\'e}ristique nulle. On donne, dans cette section,
 un proc{\'e}d{\'e}  g{\'e}n{\'e}ral qui associe {\`a} tout $K^s$-rev{\^e}tement
de courbes~$\varphi : \cC \to \cB$, une $K^s$-surface quasi-projective lisse et irr{\'e}ductible
en conservant les propri{\'e}t{\'e}s de modules et de d{\'e}finition.

\begin{theoreme}\label{theoreme:quasi}
Soit~$K$ un corps de caract{\'e}ristique nulle,~$\cB$ une $K$-courbe 
lisse projective et g{\'e}om{\'e}triquement irr{\'e}ductible
 et~$\varphi :\cC\rightarrow \cB\otimes_KK^s$ un $K^s$-rev{\^e}tement.
Il existe une $K^s$-surface quasi-projective lisse et irr{\'e}ductible,      ayant m{\^e}me
corps des modules et m{\^e}mes corps de d{\'e}finition que $\varphi$.
\end{theoreme}

\subsection{La construction de la surface} Tout d'abord, gr{\^a}ce aux
lemmes~\ref{lem:grandgenre} et~\ref{lem:sansauto}, on peut supposer que la base~$\cB\otimes_KK^s$ est de
genre au moins~$2$ et qu'elle n'a pas de $K^s$-automorphismes non triviaux.


On consid{\`e}re  $\cP=(\cB\otimes_KK^s)\times \cC$ le produit cart{\'e}sien 
des deux courbes  et  $\cG$ le graphe de $\varphi$
dans $\cP$. Soit $\cU$ l'ouvert de $\cP$ compl{\'e}mentaire
de $\cG$. On note $p_\cC : \cP\rightarrow \cC$  et $p_\cB :\ cP\rightarrow
\cB$ 
les projections sur chacun des deux facteurs.

La surface attendue n'est rien d'autre que l'ouvert~$\cU$ que l'on appelle {\em empreinte}
du rev{\^e}tement $(\cC,\varphi)$. Il nous reste {\`a} {\'e}tudier les corps des modules et de d{\'e}finition
de cette empreinte en fonction de ceux de~$\varphi$.

\subsection{Corps des modules et corps de d{\'e}finition}

On proc{\`e}de en encha{\^\i}nant deux lemmes.

\begin{lem}\label{lem:morphismeempreinte}
Soient~$\cU$ et~$\cV$ les empreintes respectives de deux $K^s$-rev{\^e}tements
connexes non
triviaux~$\varphi : \cC \to \cB\otimes_KK^s$ et~$\psi : \cD \to
\cB\otimes_KK^s$,
o{\`u}~$\cB$ est une $K$-courbe lisse, projective, g{\'e}om{\'e}triquement irr{\'e}ductible, de genre au moins~$2$
et sans $K^s$-automorphisme non trivial.
Alors tout  $K^s$-morphisme surjectif de~$\cU \to \cV$  est de la forme
$\Id \times \gamma$ o{\`u} $\gamma :\cC \rightarrow \cD$ est un  $K^s$-morphisme
entre les rev{\^e}tements~$\varphi : \cC \to \cB\otimes_KK^s$
 et~$\psi : \cD \to \cB\otimes_KK^s$.
\end{lem}

\begin{preuve}
Un $K^s$-morphisme de  $(\cC,\varphi)$ dans $(\cD,\psi)$ est un $K^s$-morphisme de courbes 
$\gamma :\cC\rightarrow\cD$  tel que 
$\psi\gamma= \varphi$. L'isomorphisme produit $\Id\times \gamma$ de 
$(\cB\otimes_KK^s)\times \cC$ dans $(\cB\otimes_KK^s)\times \cD$ 
envoie le graphe de $\varphi$ sur celui
de $\psi$ et   $\cU$ sur $\cV$.

R{\'e}ciproquement, soit $\upsilon$ un $K^s$-morphisme surjectif de $\cU$ sur $\cV$. Soit $c$ un $K^s$-point ferm{\'e}
 de $\cC$. La restriction de   $p_\cD\circ \upsilon$ {\`a}
$((\cB\otimes_KK^s)\times \{c\})\cap \cU$ est constante car le genre de $\cB$ est inf{\'e}rieur {\`a} celui de
$\cD$. On note $\gamma (c)$ cette constante et on d{\'e}finit ainsi un morphisme 
 $\gamma : \cC \rightarrow \cD$ qui n'est pas constant car $\upsilon$ est
 surjectif.
 La
 restriction de   $p_\cB\circ \upsilon$ {\`a}
$((\cB\otimes_KK^s)\times \{c\})\cap \cU$ est    un morphisme  $\beta_c$ {\`a} valeurs  dans $\cB$. Soit
$F\subset \cC$  l'ensemble des $K^s$-points ferm{\'e}s de $\cC$ tels que $\beta_c$ est constante. C'est 
un ensemble fini car $\upsilon$ est surjectif. Pour un $K^s$-point ferm{\'e} $c\not \in F$  le morphisme 
$\beta_c$  induit   un automorphisme de $\cB$, qui est
trivial puisque $\cB$ n'a pas d'automorphisme non-trivial. On a donc 
 $\upsilon (b,c)=(b,\gamma (c))$ pour  les $K^s$-points ferm{\'e}s $b$ de $\cB$ et $c$
de $\cC$ avec $c \not \in F$ et $(b,c)\in \cU$.    Soit $b$ un $K^s$-point ferm{\'e} de $\cB\otimes_KK^s$. La restriction de $p_\cB \circ \upsilon$
{\`a} $(\{b \}\times \cC)\cap \cU$ est constante {\'e}gale {\`a} $b$ sur l'ouvert non vide
$  (\{b \}\times (\cC - F))\cap \cU  $.  Donc elle est constante. Ainsi $F$ est vide et $\upsilon $
est la restriction de $\Id \times \gamma$ {\`a} $\cU$. Donc
$\Id\times \gamma$
envoie $\cU$ sur $\cV$ donc  $\psi\gamma=\varphi$.
\end{preuve}

\begin{lem}\label{lem:produitdefi}
Soit~$\cU$ l'empreinte d'un $K^s$-rev{\^e}tement~$\varphi : \cC \to \cB\otimes_KK^s$ de degr{\'e} au moins~$2$,
o{\`u}~$\cB$ est une $K$-courbe
lisse, projective, g{\'e}om{\'e}triquement irr{\'e}ductible, de genre au moins~$2$ et sans $K^s$-automorphisme
non trivial. Alors~:
\begin{enumerate}
\item le groupe des $K^s$-automorphismes de $\cU$  est {\'e}gal au groupe
des $K^s$-automorphismes du rev{\^e}tement $\varphi : \cC \to \cB\otimes_KK^s$~;
\item le corps des modules de $\cU$ (dans la cat{\'e}gorie des vari{\'e}t{\'e}s quasi-projectives)
est le corps des modules du rev{\^e}tement $\varphi : \cC \to \cB\otimes_KK^s$~;
\item une extension alg{\'e}brique de $K$ est corps de d{\'e}finition de~$\cU$ 
si et seulement si elle est corps de d{\'e}finition du rev{\^e}tement~$\varphi : \cC\to\cB\otimes_KK^s$.
\end{enumerate}
\end{lem}

\begin{preuve}
Les deux premiers points r{\'e}sultent du lemme \ref{lem:morphismeempreinte}.
On prouve le troisi{\`e}me en appliquant le principe g{\'e}n{\'e}ral expos{\'e} dans l'introduction.
Soient $L\subset K^s$ une extension alg{\'e}brique
 de $K$ et $\psi : \cD\rightarrow \cB\otimes_KL$ un $L$-mod{\`e}le
de $\varphi$. Il existe $J : \cC\rightarrow \cD\otimes_LK^s$ un  $K^s$-isomorphisme tel que $\varphi=\psi
J$. Soit
$\cV$ l'ouvert de $(\cB\otimes_KL)\times \cD$ compl{\'e}mentaire du graphe de $\psi$. C'est
une
$L$-vari{\'e}t{\'e} et $ \Id \times J$ est un  $K^s$-isomorphisme de $\cU$ sur $\cV\otimes_LK^s$.

R{\'e}ciproquement, soit $\cV$ un $L$-mod{\`e}le de $\cU$ et $\upsilon : \cU
\rightarrow
\cV\otimes_LK^s$ un $K^s$-isomorphisme. Pour tout $\sigma$ dans $\Gal(K^s/L)$
la compos{\'e}e $w_\sigma=\upsilon ^{-1}{}^\sigma \upsilon$ est un $K^s$-isomorphisme
de ${}^\sigma \cU$ sur $\cU$ qui s'{\'e}tend en un $K^s$-isomorphisme $ \Id \times
V_\sigma$ de $(\cB\otimes_K K^s)\times {}^\sigma \cC$ dans $ (\cB\otimes_KK^s)\times \cC$ avec
$V_\sigma$ un $K^s$-isomorphisme de rev{\^e}tements de ${}^\sigma \varphi$
sur $\varphi$. Puisque les
$w_\sigma$ satisfont les conditions de Weil pour la descente {\`a} $L$, il en va de m{\^e}me des $V_\sigma$.
Donc $\varphi$ admet un $L$-mod{\`e}le.
\end{preuve}

\medbreak

Le th{\'e}or{\`e}me~\ref{theoreme:catrevetement} nous permet donc d'affirmer que~:

\begin{corollaire}
Il existe des $\overline{\QQ}$-surfaces quasi-projectives lisses et g{\'e}om{\'e}triquement irr{\'e}ductibles de
corps des modules~$\QQ$, d{\'e}finies sur tous les compl{\'e}t{\'e}s de~$\QQ$ mais n'{\'e}tant pas d{\'e}finies
sur~$\QQ$.
\end{corollaire}

\section{Surfaces normales  projectives}\label{section:surfacesnormalesprojectives}

Dans cette section $K$ est encore un corps de caract{\'e}ristique nulle.
Dans le m{\^e}me esprit que la section pr{\'e}c{\'e}dente, on donne ici un proc{\'e}d{\'e} assez
g{\'e}n{\'e}ral qui associe {\`a} un
rev{\^e}tement de courbes
de corps des modules $K$, 
une $K^s$-surface  projective normale, irr{\'e}ductible 
ayant m{\^e}me corps des modules et   m{\^e}mes corps 
de d{\'e}finition que le rev{\^e}tement initial. Plus pr{\'e}cis{\'e}ment on veut montrer le~:


\begin{theoreme}\label{theoreme:projective}
Soit $\cB$ une courbe d{\'e}finie sur 
$K$ et soit $\varphi :\cC\rightarrow \cB\otimes_KK^s$ un $K^s$-rev{\^e}tement
de corps des modules~$K$.
Il existe une $K^s$-surface  projective,  irr{\'e}ductible, et normale, de corps
des modules $K$,   et 
  ayant  m{\^e}mes corps de d{\'e}finition que $\varphi$.
\end{theoreme}

\subsection{La construction de la surface~$\cS$ et ses premi{\`e}res propri{\'e}t{\'e}s}
\label{section:constructionsurfacenormale}


Notre point de d{\'e}part est encore un 
rev{\^e}tement~$\varphi : \cC \to \cB\otimes_K L$ de degr{\'e}~$d\ge 2$, de corps des
modules~$K$ et d{\'e}fini sur~$L\subset K^s$ une extension alg{\'e}brique
 de $K$.
D'apr{\`e}s le lemme~\ref{lem:sansauto}, on peut supposer que $\cB$ est de genre au moins~$2$ et
qu'elle n'a pas de $\Qb$-automorphisme non trivial.
On construit un rev{\^e}tement de $(\cB\otimes_KL)\times \cD$
fortement ramifi{\'e} le long du graphe de $\varphi$.

On se donne  un syst{\`e}me fini $(h_j)_{1\le j\le J}$ de g{\'e}n{\'e}rateurs
du corps $K(\cB)$. On suppose qu'aucun  $h_j$ n'est une puissance dans 
$K^s(\cB)$.

On pose~$I = 2J$ et~$\Pi = \prod_{i = 1}^{I} p_i$  le  produit
des $I$ premiers nombres premiers plus grands que le degr{\'e} $d$ de $\varphi$.

On choisit deux entiers positifs
$a$ et $b$ premiers entre eux et plus grands que 
$1+2(g(\cB)-1)\Pi+\Pi^2$ et pour tout $1\le j\le J$ on pose~:
$$
f_j=h_j^a,
\qquad
f_{j+J}=h_j^b.
$$
Quitte {\`a} grandir les valeurs de~$a$ et~$b$, on peut supposer qu'aucune des
fonctions~$f_i - \lambda$ n'est
une puissance $p_i$-{\`e}me dans $K^s(\cB)$ 
pour~$\lambda \in K^s$ et~$1 \leq i \leq I$ :
c'est {\'e}vident pour $\lambda=0$. Si $\lambda\not =0$
et  si $h_i^a-\lambda=\prod_{0\le k \le a-1}
(h_i-\zeta_{a}^k\lambda^{\frac{1}{a}})$ est une puissance,  alors 
$h_i$ a au moins  $a$ valeurs singuli{\`e}res distinctes. Ceci est
exclu si $a$ est plus grand que le nombre de valeurs singuli{\`e}res de $h_i$.

 D'autre part, les $(f_i)_{1\le i\le I}$ engendrent 
$K(\cB)$  et elles  ont un 
degr{\'e} plus grand que $1+2(g(\cB)-1)\Pi+\Pi^2$.
Soit~$M$ le maximum et $m$ le minimum des degr{\'e}s des $f_i$, alors~:
$$
\forall 1 \leq i \leq I, \qquad 1+2(g(\cB)-1)\Pi+\Pi^2 < m \leq \deg(f_i) \leq M.
$$

Soit $p$ un nombre premier plus grand que $(g(\cB)+IM)\Pi$ et
notons $\cD$ la $L$-courbe et $\psi : \cD\rightarrow \cB\otimes_KL$
le rev{\^e}tement  de degr{\'e} $pd$, donn{\'e}s par le lemme \ref{lem:granddegre}.
Le genre de $\cD$ est au moins $dp>(g(\cB)+IM)\Pi$ et les rev{\^e}tements~$\varphi$ et~$\psi$ ont m{\^e}me
corps des modules et m{\^e}mes corps de d{\'e}finition. On v{\'e}rifie en outre
qu'aucune des
fonctions~$f_i\circ \psi - \lambda$ n'est
une puissance $p_i$-{\`e}me dans $K^s(\cD)$
pour~$\lambda \in K^s$ et~$1 \leq i \leq I$ (ceci est d{\'e}j{\`a} vrai pour $f_i-\lambda$ et de plus le degr{\'e} $pd$ 
de $\psi$ est premier {\`a} $p_i$.)

On d{\'e}finit la fonction $g_i$ sur $(\cB\otimes_K L)\times \cD$
par~:
$$
g_i (P,Q)=f_i(\psi(Q))-f_i(P).
$$
La partie n{\'e}gative $(g_i)_\infty$ du diviseur de $g_i$ est~:
\begin{equation*}
(g_i)_\infty=(f_i)_\infty\times \cD+\cB\times (f_i\circ\psi)_\infty,
\end{equation*}
et les parties positives  $(g_i)_0$ que l'on note 
$\Delta_i$ v{\'e}rifient quant-{\`a}-elles~:
\begin{equation*}
\pgcd _i(\Delta_i)=\cG,
\end{equation*}
o{\`u} $\cG$ est le graphe de $\psi$.

On d{\'e}finit une extension ab{\'e}lienne du corps des fonctions 
$L((\cB\otimes_KL)\times\cD)$
de $(\cB\otimes_KL)\times\cD$ en posant 
$$
y_i^{p_i}=g_i
$$
o{\`u} $p_i$ est $i$-{\`e}me nombre premier plus grand que le degr{\'e}
$d$ de $\varphi$.

On note, enfin, $\cS$ la normalisation de $(\cB\otimes_K L)\times \cD$ dans le corps ainsi
d{\'e}fini. Par construction~$\cS$ est normale. De plus, on dispose d'un rev{\^e}tement de
degr{\'e}~$\Pi$, ramifi{\'e}~:
$$
\chi : \cS\rightarrow (\cB\otimes_KL)\times \cD.
$$

Afin de pr{\'e}parer l'{\'e}tude du groupe d'automorphismes
de $\cS$, on introduit une famille de courbes
sur $\cS$. Pour~$Q$  un point  de $\cD$,
notons~$\cE_Q$ l'image r{\'e}ciproque de $\cB \times Q$ par $\chi$
et~$\chi_Q : \cE_Q\rightarrow \cB \times Q$ le rev{\^e}tement
obtenu par restriction  de $\chi$. 
Si $Q$ est un point
 g{\'e}n{\'e}rique\footnote{Par {\flqq point g{\'e}n{\'e}rique\frqq}
de $\cD\otimes_LK^s$, nous entendons, ici et dans la suite de ce texte, un point d{\'e}fini
sur le corps des fonctions rationnelles~$K^s(\cD)$.} de $\cD\otimes_LK^s$, la courbe $\cE_Q$ est irr{\'e}ductible, et le rev{\^e}tement~$\chi_Q$
est de degr{\'e} $\Pi$. De plus, par construction, le degr{\'e}
du diviseur de ramification de ce rev{\^e}tement est  major{\'e} par le produit $2IM$
(o{\`u}, on le rappelle, $I$ est le nombre de fonctions dans la famille $(f_i)_i$
et $M$ est le maximum des degr{\'e}s de ces fonctions).
Le genre g{\'e}om{\'e}trique de $\cE_Q$ est donc major{\'e} par~:
\begin{equation} \label{majoration_genre} 
g(\cE_Q) \leq( g(\cB)+IM)\Pi < g(\cD).
\end{equation}
On dispose aussi d'une minoration du genre de n'importe quel sous-rev{\^e}tement non trivial
de~$\chi_Q$. En effet, encore par construction, un tel sous-rev{\^e}tement est de degr{\'e} au moins {\'e}gal
{\`a}~$p_1 \geq 3$ avec un diviseur de ramification de degr{\'e} au moins {\'e}gal {\`a}~$m$
(o{\`u}~$m$ est le minimum des degr{\'e}s des fonctions~$f_i$), si bien que~:
\begin{equation} \label{minoration_genre} 
1+2(g(\cB)-1)\Pi+\Pi^2 < m \leq g(\text{sous-rev{\^e}tement de~$\chi_Q : \cE_Q \to \cB$}).
\end{equation}

Les in{\'e}galit{\'e}s (\ref{majoration_genre}) et 
(\ref{minoration_genre}) seront utiles pour calculer le groupe d'automorphismes de $\cS$.

\subsection{\'Etude du groupe d'automorphismes de~$\cS$}

Notons~$A$ le groupe des $K^s$-automorphismes de~$\psi$. Vus comme des $K^s$-automorphismes de la
surface~$(\cB\otimes_KL)\times\cD$, ces automorphismes se rel{\`e}vent en des automorphismes de $K^s(\cS)/K^s$ qui
fixent les $y_i$ et qui stabilisent $K^s( \cB\otimes_KL\times \cD)$
(dans la suite, nous adopterons la m{\^e}me notation pour un automorphisme de~$\psi$ et ses rel{\`e}vements
{\`a}~$(\cB \otimes_KL)\times \cD$ et {\`a}~$\cS$).
Autrement dit, tout {\'e}l{\'e}ment~$\alpha \in A$ induit un automorphisme  de~$\cS$ et on peut
voir~$A$ comme un sous-groupe de~$\Aut_{K^s}(\cS)$, le groupe de $K^s$-automorphismes de~$\mathcal{S}$.
On dispose d'un autre sous-groupe de~$\Aut_{K^s}(\cS)$, {\`a} savoir~$B=\prod_i\ZZ/p_i\ZZ$ le groupe de
Galois de l'extension $K^s(\cS)/K^s((\cB\otimes_KL) \times\cD)$.

En r{\'e}sum{\'e}, les {\'e}l{\'e}ments~$\alpha \in A$,  vus
 comme des {\'e}l{\'e}ments de~$\Aut_{K^s}(\cS)$
font commuter le diagramme suivant~:
$$
\xymatrix@=15pt{
\cS\otimes_LK^s \ar@{->}[r]^{\alpha} \ar@{->}[d]_{\chi} & \cS \otimes_LK^s\ar@{->}[d]^{\chi} \\
((\cB\otimes_KL) \times \cD)\otimes_LK^s \ar@{->}[r]^{\alpha} & 
((\cB\otimes_KL) \times \cD)\otimes_LK^s
}
$$

Les {\'e}l{\'e}ments~$\beta \in B$,  vus
 comme des {\'e}l{\'e}ments de~$\Aut_{K^s}(\cS)$
font commuter le diagramme suivant~:
$$
\xymatrix@=15pt{
\cS\otimes_LK^s \ar@{->}[rr]^{\beta} \ar@{->}[rd]_{\chi} & & \cS 
\otimes_LK^s\ar@{->}[ld]^{\chi} \\
& ((\cB\otimes_KL) \times \cD)\otimes_LK^s
}
$$

Et il est clair que~$A \times B \subset \Aut_{K^s}(\cS)$. En fait l'inclusion inverse est v{\'e}rif{\'e}e.

\begin{lem}\label{lem:automorphismes}
Le groupe des $K^s$-automorphismes
de $\cS$ est {\'e}gal {\`a} $A\times B$.
\end{lem}

\begin{preuve}
Soit~$\agot$ un $K^s$-automorphisme  de $\cS$.

Tout d'abord, si~$Q$ est un point g{\'e}n{\'e}rique de~$\cD \otimes_L K^s$, on sait gr{\^a}ce
{\`a} l'in{\'e}galit{\'e}~(\ref{majoration_genre}) de la section~\ref{section:constructionsurfacenormale},
que~$g(\cE_Q) < g(\cD)$.
Il en r{\'e}sulte que l'automorphisme~$\agot$ est tel que~$\agot(\cE_Q)=\cE_{\alpha(Q)}$
o{\`u} $\alpha$ est un $K^s$-automorphisme de $\cD$.


V{\'e}rifions  que l'isomorphisme de~$\cE_Q \to \cE_{\alpha(Q)}$ induit par~$\agot$
fait commuter le diagramme~:
$$
\xymatrix@=30pt{
\cE_Q \ar@{->}[r]^{\agot} \ar@{->}[d]_{\chi_Q} & \cE_{\alpha(Q)} \ar@{->}[d]^{\chi_{\alpha(Q)}} \\
\cB \times Q \ar@{->}[r]^{ id \times \alpha} & \cB \times \alpha(Q)
}
$$
avec~${\rm id} \times \alpha$ comme isomorphisme en bas.
Le produit des  fonctions $\chi_Q$ et $\chi_{\alpha (Q)}\circ \agot$ 
d{\'e}finit un morphisme~:
$$
\cE_Q 
\xrightarrow{\chi_Q \times \chi_{\alpha (Q)}\circ \agot} \cB \times \cB,
$$
dont l'image
est un diviseur de bidegr{\'e} $\le (\Pi,\Pi)$. D'apr{\`e}s le lemme \ref{lem:adjonction},
le genre arithm{\'e}tique de cette image est au plus $1+2(g(\cB)-1)\Pi+\Pi^2$ . D'autre part, en
composant avec la premi{\`e}re projection de~$\cB \times \cB \to \cB$,
on s'aper{\c c}oit que cette image est {\it coinc{\'e}e} entre $\cE_Q$
et $\cB\otimes_KK^s$. Mais on a vu, dans l'in{\'e}galit{\'e}~(\ref{minoration_genre}) de la
section~\ref{section:constructionsurfacenormale}, qu'un tel sous-rev{\^e}tement non-trivial
de $\chi_Q :\cE_Q\rightarrow \cB\otimes_KK^s$ a un genre g{\'e}om{\'e}trique
au moins $m>1+2(g(\cB)-1)\Pi+\Pi^2$. Il en r{\'e}sulte que l'image  de 
$\chi_Q\times (\chi_{\alpha (Q)}\circ \agot)$ est forc{\'e}ment $K^s$-isomorphe
{\`a} $\cB$ si bien que c'est une correspondance de bidegr{\'e} $(1,1)$. Comme $\cB$ n'a
pas de $K^s$-automorphisme non-trivial on en d{\'e}duit que le diagramme pr{\'e}c{\'e}dent peut {\^e}tre
compl{\'e}t{\'e} en bas avec l'identit{\'e}. Donc $\chi\circ \agot=(\Id\times
\alpha)\circ
\chi$.

Montrons maintenant  que~$\alpha \in A$. Pour cela, remarquons que pour~$Q$
un point g{\'e}n{\'e}rique de $\cD\otimes_KK^s$, on vient  de montrer que~$\agot$ induit un isomorphisme entre les
rev{\^e}tements~$\chi_Q : \cE_Q \to \cB$ et~$\chi_{\alpha(Q)} : \cE_{\alpha(Q)} \to \cB$.
Ces deux rev{\^e}tements ont donc m{\^e}mes donn{\'e}es  de ramification, c'est-{\`a}-dire respectivement les
points~$P$ tels que~$f_i(P) = f_i(\psi(Q))$ et~$f_i(P) = f_i(\psi(\alpha(Q)))$.
Ainsi~:
$$
\forall i, \quad f_i(\psi(Q)) = f_i(\psi(\alpha(Q)))
$$
\noindent donc $\psi(Q) = \psi(\alpha(Q))$, 
car les~$f_i$ 
engendrent le corps $K(\cB)$. Autrement dit~$\psi = \psi \circ \alpha$, ce
qui signifie que~$\alpha\in A$. Donc $\chi\circ \alpha=\alpha\circ
\chi=\chi\circ \agot$. 
On pose   $\beta= \agot\circ \alpha^{-1}$ et on v{\'e}rifie que $\chi\circ
\beta=\chi$ donc $\beta$ est dans $B$. 
En conclusion, on a bien~$\agot=\beta \alpha \in A\times B$.
\end{preuve}

\begin{remarque}
La preuve ci-dessus montre que le groupe des $K^s$-automorphisme
birationnels de $\cS$ est {\'e}gal {\`a} $A\times B$. Nous n'aurons pas besoin
de ce r{\'e}sultat.
\end{remarque}

\subsection{Corps des modules et corps de d{\'e}finition de~$\cS$} 
Tout d'abord, le corps des
modules de $\cS$ est $K$ car un isomorphisme entre ${}^\sigma \varphi$ et 
$\varphi$ donne lieu {\`a} un isomorphisme entre ${}^\sigma\cS$ et $\cS$.

Ensuite, la surface~$\cS$ {\'e}tant construite {\`a} partir de $\varphi : \cC \rightarrow \cB\otimes_KL$ 
sans extension des scalaires, tout corps de d{\'e}finition de~$\varphi$ est corps de
d{\'e}finition de~$\cS$. 

Il nous reste {\`a} v{\'e}rifier que tout corps
de d{\'e}finition de~$\cS$ est aussi
un corps de d{\'e}finition de~$\varphi$ ou, ce qui revient au m{\^e}me, de~$\psi$.
Pour cela, consid{\'e}rons ${}^\sigma\cS$ et  ${}^\tau\cS$  deux
conjugu{\'e}s de  $\cS$
par action de Galois et $I : {}^\tau\cS \rightarrow {}^\sigma\cS  $  un
$K^s$-isomorphisme.
Comme le groupe de $K^s$-automorphismes de la surface $\cS$  admet un unique sous-groupe
d'ordre $\Pi$, il en va de m{\^e}me pour chacune de ses deux conjugu{\'e}es et les deux
groupes
en questions se correspondent par $I$. On obtient donc par quotient, un
$K^s$-isomorphismes de ${}^\tau((\cB\otimes_K L)\times \cD)$  dans 
${}^\sigma((\cB\otimes_K L)\times \cD)$ qui envoie le graphe de~${}^\tau \psi$ dans
celui de ${}^\sigma \psi$ (donn{\'e}es de ramification
des rev{\^e}tements quotients) et donc l'empreinte de~${}^\tau \psi$ dans celle
de~${}^\sigma \psi$. Ce proc{\'e}d{\'e} permet de d{\'e}finir un foncteur
covariant
qui transforme
une donn{\'e}e de descente 
pour~$\cS$ en une donn{\'e}e de descente
pour l'empreinte du rev{\^e}tement~$\psi$ sur la
surface~$(\cB \otimes_K L) \times \cD$.
Gr{\^a}ce au th{\'e}or{\`e}me~\ref{theoreme:criterech}, cela montre que 
tout 
corps de d{\'e}finition
de $\cS$,  est aussi 
corps de d{\'e}finition de l'empreinte du rev{\^e}tement~$\psi$. En vertu du
lemme~\ref{lem:produitdefi}, on sait qu'alors
ce corps  est aussi  corps de d{\'e}finition  de~$\psi$. 

\medbreak

Gr{\^a}ce, une nouvelle fois, au th{\'e}or{\`e}me~\ref{theoreme:catrevetement}, nous en d{\'e}duisons que~:

\begin{corollaire}
Il existe des $\overline{\QQ}$-surfaces normales, projectives et g{\'e}om{\'e}triquement irr{\'e}ductibles, de corps des modules
{\'e}gal {\`a}~$\QQ$, d{\'e}finies sur tous les compl{\'e}t{\'e}s de~$\QQ$, mais ne pouvant pas {\^e}tre d{\'e}finies
sur~$\QQ$.
\end{corollaire}

\section{Courbes}\label{section:courbes}

Dans cette section $K$ est encore un corps de caract{\'e}ristique nulle.
Dans le m{\^e}me esprit que la section pr{\'e}c{\'e}dente, on donne ici un proc{\'e}d{\'e} assez
g{\'e}n{\'e}ral qui associe {\`a} un
rev{\^e}tement de courbes
de corps des modules $K$, 
une $K^s$-courbe  projective lisse, irr{\'e}ductible,
de corps des modules $K$, ayant les  m{\^e}mes corps 
de d{\'e}finition que le rev{\^e}tement initial. Plus pr{\'e}cis{\'e}ment on veut montrer le~:


\begin{theoreme}\label{theoreme:courbe}
Soit $\cB$ une courbe d{\'e}finie sur 
$K$ et soit $\varphi :\cC\rightarrow \cB\otimes_KK^s$ un $K^s$-rev{\^e}tement
de corps des modules~$K$.
Il existe une $K^s$-courbe  projective,  irr{\'e}ductible, et lisse, de corps
des modules $K$ et 
  ayant  m{\^e}mes corps de d{\'e}finition que $\varphi$.
\end{theoreme}

Le point de d{\'e}part est la surface~$\mathcal{S}$ construite {\`a} la
section \ref{section:surfacesnormalesprojectives}. On garde donc les notations
de cette pr{\'e}c{\'e}dente section. On sait que $\cS$
 a pour corps des modules~$K$ et pour corps
de d{\'e}finition les corps de d{\'e}finition du rev{\^e}tement 
$\varphi : \cC/L \rightarrow \cB\otimes _KL$
(ou $\psi : \cD/L \rightarrow \cB\otimes _KL$,  cela revient au m{\^e}me).

L'id{\'e}e est de tracer sur~$\mathcal{S}$ une courbe singuli{\`e}re (mais stable) 
h{\'e}ritant des propri{\'e}t{\'e}s de module et de d{\'e}finition
de~$\mathcal{S}$ puis de d{\'e}former cette courbe pour
aboutir {\`a} une courbe projective lisse tout en conservant {\'e}videmment les propri{\'e}t{\'e}s de module et de
d{\'e}finition.

\subsection{Deux courbes stables}\label{subsection:deux}

Les notations sont celles de la section~\ref{section:constructionsurfacenormale}. Soit~$\Gamma$
l'union de tous les supports des diviseurs des fonctions~$g_i$. Elle contient le lieu de
ramification du
rev{\^e}tement~$\chi$.


\begin{lem}
Il existe deux fonctions non constantes~$f,g \in K(\mathcal{B})$ 
telles que~:

\begin{itemize}
\item[$\bullet$] le diviseur~$((f)_0 + (f)_\infty) \times \mathcal{D}$ coupe transversalement~$\Gamma$~;
\item[$\bullet$] le diviseur~$(\mathcal{B}\otimes_KL) \times ((g \circ
  \psi)_0+(g\circ \psi)_\infty)$ coupe
transversalement~$\Gamma \cup [((f)_0 +(f)_\infty )\times \cD]$;
\item[$\bullet$] tout $K^s$-automorphisme de~$\cD$ qui stabilise la fibre~$(g \circ \psi)_0$
est en fait un automorphisme du rev{\^e}tement~$\psi$ (notons que cette fibre
est simple en vertu de la condition pr{\'e}c{\'e}dente);
\item[$\bullet$] Pour tout  z{\'e}ro $P$ de $f$, on note $\chi_P : 
\cF_P \rightarrow P\times \cD$   la  correstriction\footnote{En toute rigueur
il faudrait {\'e}crire $P\times (\cD\otimes_LK^s)$.}
de $\chi$ {\`a} $P\times \cD$ et on demande que le rev{\^e}tement
$\kappa_P:= g\circ \psi \circ \chi_P  : \cF_P \rightarrow \PU$ n'ait
 pas d'autres
automorphismes que les {\'e}l{\'e}ments de $A\times B$ :

$$\Aut_{K^s}(g\circ \psi \circ \chi_P)=\Aut_{K^s}(\psi\circ \chi_P)=A\times B.$$
\end{itemize}
\end{lem}

\begin{preuve}
Soit $f\in K (\cB)$ une fonction non-constante. 
On applique  le lemme \ref{lem:fibrescoupenttransversalement} 
{\`a} $L$, $\cB$,   $\cD$, $\Gamma$ et $f$. On en d{\'e}duit qu'il existe
un $x$ et un $y$ distincts dans $K$ tels que $(f)_x\times \cD$ et $(f)_y\times
\cD$ coupent  transversalement $\Gamma$.
On remplace $f$ par $(f-x)/(f-y)$ et la premi{\`e}re condition est satisfaite.

On observe  alors  que pour tout z{\'e}ro $P$  de $f$,  la courbe $\cF_P$ est lisse
et g{\'e}om{\'e}triquement irr{\'e}ductible car $(f)_0\times \cD$ coupe transversalement
le lieu $\Gamma$ de ramification  de $\chi$. En outre 

$$\Aut_{K^s}(\psi\circ\chi_P)=A\times B$$
\noindent  car $K^s(\cF_P)/K^s(\cB)$ est le compositum de
$K^s(\cD)/K^s(\cB)$
et de l'extension ab{\'e}lienne  de $K^s(\cB)$ engendr{\'e}e par les
$g_i(P,)^{\frac{1}{p_i}}= (f_i\circ \psi-f_i(P))^{\frac{1}{p_i}}$. 

On cherche une fonction $g$ dans $K(\cB)$ telle que $g\circ \psi$
n'ait  pas d'autre $K^s$-automorphisme  que les {\'e}l{\'e}ments de $A$
et, pour tout z{\'e}ro  $P$ de $f$, le rev{\^e}tement $\kappa_P=
g\circ \psi \circ \chi_P$ n'ait pas d'autre $K^s$-automorphisme  
que les {\'e}l{\'e}ments de $\Aut_{K^s}(\psi\circ\chi_P)=A\times B$. 
D'apr{\`e}s le lemme 
\ref{lem:fonctionsauto}, les fonctions de $K^s(\cB)$ qui ne satisfont
pas toutes ces hypoth{\`e}ses sont contenues dans une union finie de
sous-$K$-alg{\`e}bres  strictes.
Il existe donc une telle fonction $g$. 

D'apr{\`e}s le lemme
\ref{lem:fibrescoupenttransversalement},
les $x$ dans $K$ tels que $(g\circ\psi)_x$ ne coupe pas 
$\Gamma \cup [((f)_0 +(f)_\infty ) \times \cD]$ transversalement sont en nombre fini.

D'apr{\`e}s le lemme
\ref{lem:autofibres},
les $x$ dans $K$ tels que $(g\circ\psi)_x$ ait
un $K^s$-automorphisme en dehors de $A$ 
sont en nombre fini.

Donc il existe   $x$ et $y$ distincts dans $K$ tels que 
$(g\circ\psi)_x$  et $(g\circ\psi)_y$  coupent
$\Gamma \cup [((f)_0 + (f)_\infty )\times \cD]$ transversalement et
$(g\circ\psi)_x$ n'ait
pas d'autre  automorphisme que les {\'e}l{\'e}ments de  $A$.
On remplace $g$ par $(g-x)/(g-y)$ et les trois derni{\`e}res conditions sont
satisfaites.
\end{preuve}

\smallbreak

Soit alors~$\mathcal{J}_0$ la courbe, trac{\'e}e sur~$(\mathcal{B} \otimes_K L) \times \mathcal{D}$, d'{\'e}quation~:
$$
f(P) \times g\circ \psi(Q) = 0.
$$
On peut montrer le r{\'e}sultat suivant.

\begin{lem}\label{lem:automorphismes2}
La courbe~$\mathcal{J}_0$ est stable et telle que~$\Aut_{K^s}(\mathcal{J}_0) \simeq A$.
\end{lem}

\begin{preuve}
La courbe~$\mathcal{J}_0$ est r{\'e}duite car les z{\'e}ros de~$f$ et~$g\circ \psi$ sont simples
et ses points singuliers ---~c'est-{\`a}-dire les~$(P,Q)$ sur
$(\mathcal{B} \otimes_K K^s) \times \mathcal{D}$ tels que~$f(P) = g\circ \psi(Q) = 0$~--- sont des points doubles ordinaires~;
elle est donc semi-stable. Elle est aussi connexe et projective. Comme toutes ses composantes irr{\'e}ductibles sont
isomorphes {\`a}~$\mathcal{B}$ ou~$\mathcal{D}$ toutes deux de genre~$\geq 2$, la
courbe~$\mathcal{J}_0$ est bien une courbe stable.

Enfin~$\Aut_{K^s} (\mathcal{J}_0) \simeq A$~: le groupe des $K^s$-automorphismes 
de $\cJ_0$ est le groupe $A$ des automorphismes
 de $\psi$.
En effet, un automorphisme $\agot$ de $\cJ_0$ permute les composantes
 irr{\'e}ductibles. Certaines de ces composantes sont des copies de $\cB$ et
 d'autres sont des copies de $\cD$. Comme $\cB$ et $\cD$ ne sont pas
 isomorphes, $\agot$ respecte les deux sous-ensembles.  La restriction de
 $\agot$
{\`a} une composante isomorphe {\`a} $\cB$, suivie de la projection sur $\cB$ est un
 morphisme non-constant de $\cB$ dans elle-m{\^e}me, donc l'identit{\'e} car $\cB$ n'a
pas d'automorphismes. Donc $\agot$ stabilise chaque composante isomorphe {\`a}
 $\cD$. Les points singuliers sur une telle composante sont les z{\'e}ros de
 $g\circ \psi$. 
Comme l'ensemble de ces z{\'e}ros n'est
pas stabilis{\'e} par d'autres automorphismes 
 que ceux
de $\psi$, on peut supposer que $\agot$  est trivial sur  une
des composantes isomorphes {\`a}
 $\cD$ (quitte {\`a} le composer avec un {\'e}l{\'e}ment de $A$). Il doit donc 
 stabiliser les composantes isomorphes {\`a} $\cB$. Et comme elles n'ont pas
 d'automorphisme, il agit trivialement dessus. Comme $\agot$ agit fid{\`e}lement
sur les z{\'e}ros de $g\circ \psi$, il agit trivialement
sur toutes  les composantes isomorphes
{\`a} $\cD$.
\end{preuve}

Soit~$\mathcal{K}_0$ l'image r{\'e}ciproque de~$\mathcal{J}_0$ par~$\chi$. On
{\'e}tudie {\`a} nouveau sa stabilit{\'e} et on cerne un sous-groupe
de son groupe d'automorphismes : le sous-groupe des automorphismes 
{\it d{\'e}formables},
not{\'e}~$\Aut_{K^s}^{\text{d{\'e}f.}} (\mathcal{K}_0)$. 
Nous expliquons maintenant le sens de ce terme.

Remarquons d'abord que les 
composantes de $\cK_0\otimes_LK^s$ sont de deux sortes. Certaines sont des rev{\^e}tements 
de $(\cB\otimes_KK^s)\times Q$ pour $Q$ un $K^s$-z{\'e}ro de $g\circ\psi$ et on les note $\cE_Q$. Les
autres sont des rev{\^e}tements 
de $P\times (\cD\otimes_KK^s)$ pour $P$ un $K^s$-z{\'e}ro de $f$ et on les
note $\cF_P$. On note $\chi_P : \cF_P\rightarrow P\times \cD$ et
$\chi_Q : \cE_Q \rightarrow \cB\times Q$ les restrictions de $\chi$
aux composantes de $\cK_0$. Soit   $T$  un point singulier de 
$\cK_0$ tel que $\chi(T)=(P,Q)$. Donc $T$ est dans l'intersection
de $\cE_Q$ et $\cF_P$. Le point de $\cE_Q$ correspondant {\`a} $T$
est appel{\'e} $R$. Le point de $\cF_P$ correspondant {\`a} $T$
est appel{\'e} $S$.
Donc $\chi_Q(R)=P$ et $\chi_P(S)=Q$.
Donc $f\circ \chi_Q$ est une uniformisante
pour $\cE_Q$ en $R$ et $g\circ \psi \circ \chi_P$ est une uniformisante
pour $\cF_P$ en $S$. Si $\agot$ est un automorphisme de $\cK_0$ et
$T'=(R',S')$ l'image de  $T=(R,S)$ par $\agot$ on note $\chi(R')=(P',Q')$.
Notons que $f\circ \chi_{Q'}\circ \agot$ est une uniformisante
pour $\cE_Q$ en $R$ et $g\circ \psi \circ \chi_{P'}\circ \agot$
 est une uniformisante
pour $\cF_P$ en $S$.

On appelle
{\it automorphisme d{\'e}formable}
de $\cK_0$
un automorphisme $\agot$ tel qu'en tout  point singulier $T$
 de 
$\cK_0$ on ait

\begin{equation}\label{eq:deformable}
\frac{f\circ \chi_{Q'}\circ \agot}{f\circ \chi_Q}(R)\times 
\frac{g\circ \psi\circ \chi_{P'} \circ \agot}{g\circ \psi  \circ \chi_P}(S)=1
\end{equation}
\noindent o{\`u} $\chi (T)=(P,Q)$ et $\chi(\agot(T))=(P',Q')$.

La justification de cette d{\'e}finition  est donn{\'e}e au paragraphe 
\ref{subsection:deformation}. Les automorphismes d{\'e}for\-mables forment
clairement un sous-groupe du groupe des automorphismes de $\cK_0$.

\begin{lem}
La courbe~$\mathcal{K}_0$ est stable et telle que~$\Aut_{K^s}^{\text{d{\'e}f.}}(\mathcal{K}_0) \simeq A \times B$.
\end{lem}

\begin{preuve}
Il est {\'e}vident que 
$A\times B$ agit (par restriction) sur 
$\cK_0$ et que les automorphismes correspondants
sont d{\'e}formables. 

La courbe~$\mathcal{K}_0$ est trac{\'e}e
sur~$\mathcal{S}$. Les conditions impos{\'e}es sur les fonctions~$f$ et~$g$ permettent de montrer
que~$\mathcal{K}_0$ est une courbe stable. En raison de la simplicit{\'e} de  l'intersection de $(f)$ et $(g)$ 
avec le lieu de ramification~$\Gamma$, les points singuliers
de~$\mathcal{J}_0$ se rel{\`e}vent en des points singuliers de~$\mathcal{K}_0$ qui restent de premi{\`e}re
esp{\`e}ce (ils sont simplement {\flqq d{\'e}multipli{\'e}s\frqq}). La connexit{\'e} de~$\mathcal{K}_0$ provient
du fait que les fonctions~$g_i$, dont on extrait une racine, ne peuvent pas {\^e}tre une puissance
$p_i$-{\`e}me car aucune des fonctions~$f_i - \lambda$, $\lambda \in K^s$, n'en
est une, et pas davantage les $f_i\circ \psi -\lambda$
(cf. \S~\ref{section:constructionsurfacenormale}).

Enfin, montrons que~$\Aut_{K^s}^{\text{d{\'e}f.}}(\mathcal{K}_0) \simeq A \times B$.
Pour le voir, remarquons que les 
composantes de $\cK_0\otimes_KK^s$ sont de deux sortes. Certaines sont des rev{\^e}tements 
de $(\cB\otimes_KK^s)\times Q$ pour $Q$ un $K^s$-z{\'e}ro de $g\circ\psi$ et
on les a not{\'e}s $\cE_Q$. Les
autres sont des rev{\^e}tements 
de $P\times (\cD\otimes_LK^s)$ pour $P$ un $K^s$-z{\'e}ro de $f$ et on les a
not{\'e}s $\cF_P$. Les $\cF_P$
et les $\cE_Q$ n'ont pas le m{\^e}me genre, donc elles ne peuvent {\^e}tre
isomorphes. Ainsi tout $K^s$-automorphisme $\agot $ de $\cK_0$ stabilise l'ensemble
des
$\cF_P$ et celui des $\cE_Q$.

Soit~$Q$ et~$Q'$ des $K^s$-z{\'e}ros de~$g \circ \psi$ tels
que~$\agot(\cE_Q) = \cE_{Q'}$. Comme dans la preuve du
lemme~\ref{lem:automorphismes}, on remarque que l'image de~$\cE_Q$ dans le
produit~$\cB \times \cB$, via le morphisme~$\chi_Q \times \chi_{Q'} \circ \agot$, a un genre
arithm{\'e}tique moindre que~$1+2(g(\cB)-1)\Pi+\Pi^2$. Cela montre encore que cette image est
 $K^s$-isomorphisme {\`a}~$\cB$ (sinon l'image rev{\^e}t~$\cB$ par un rev{\^e}tement non trivial
de genre forc{\'e}ment plus grand que~$1+2(g(\cB)-1)\Pi+\Pi^2$). Comme~$\cB$ n'a toujours pas
d'automorphisme, on en d{\'e}duit {\`a} nouveau que~$\agot$ induit un isomorphisme de rev{\^e}tements
entre les restrictions~$\chi_Q : \cE_Q \to \cB$ et~$\chi_{Q'} : \cE_{Q'} \to \cB$ de~$\chi$.
Ainsi

\begin{equation}
\chi_Q = \chi_{Q'} \circ \agot.\label{eq:chiq}
\end{equation}

D{\`e}s lors~$\agot$ stabilise chaque composante~$\cF_P$ o{\`u}~$P$ est un $K^s$-z{\'e}ro de~$f$. En effet,
partons d'un point singulier~$T =(R,S)\in  \cE_Q \cap \cF_P $, o{\`u}~$P$ est un $K^s$-z{\'e}ro de~$f$
et o{\`u}~$Q$ un $K^s$-z{\'e}ro de~$g \circ \psi$. On a donc~$\chi(T) = (P,Q) \in \cB \times \cD$.
On sait qu'il existe~$P' \in \cB(K^s)$ et~$Q' \in \cD(K^s)$ tels
que~$\agot(T) \in \cF_{P'} \cap \cE_{Q'}$ si bien que~$\agot(T) \in
\cE_{Q'}\cap \agot(\cE_{Q})$
d'o{\`u}~$\agot(\cE_{Q}) = \cE_{Q'}$. On est donc dans le contexte du paragraphe pr{\'e}c{\'e}dent,
c'est-{\`a}-dire que~:
$$
P' = \chi_{Q'} \circ \agot(T) = \chi_Q(T) = P. 
$$
Conclusion~$P = P'$ puis~$\agot(\cF_P) = \cF_P$. De plus, en vertu des 
formules~(\ref{eq:deformable}) et (\ref{eq:chiq}) on conna{\^\i}t
le quotient~:
\begin{equation} \label{gpsi_Q'_sur_gpsi_Q}
\frac{g \circ \psi \circ \chi_P \circ \agot}{g \circ \psi \circ \chi_P}(S) = 1.
\end{equation}

Notons~$\agot_P$ la restriction de~$\agot$ {\`a}~$\cF_P$. C'est donc un automorphisme de~$\cF_P$~;
mieux, montrons que c'est la restriction {\`a}~$\cF_P$ d'un {\'e}l{\'e}ment de~$A \times B$. Pour cela,
introduisons la fonction~$\kappa_P = g \circ \psi \circ \chi_P \in K^s(\cF_P)$. Elle est de
degr{\'e}~$\deg(g) \times pd \times \Pi$ et ses z{\'e}ros, tous simples, ne sont rien d'autre que les
points d'intersection de~$\cF_P$ avec les autres composantes de~$\cK_0$. Comme~$\agot_P$ permute
ces z{\'e}ros, la fonction~$\kappa_P \circ \agot_P$ a les m{\^e}mes z{\'e}ros que~$\kappa_P$ avec la m{\^e}me
multiplicit{\'e}, i.e.~$1$. D{\`e}s lors, seuls les p{\^o}les
de la fonction~$\kappa_P$ peuvent {\^e}tre p{\^o}le de la
fonction~$\frac{\kappa_P}{\kappa_P \circ \agot_P} - 1$. Cette derni{\`e}re fonction est donc de
degr{\'e} inf{\'e}rieur ou {\'e}gal
{\`a} celui de~$\kappa_P$~.
Or 
d'apr{\`e}s~$(\ref{gpsi_Q'_sur_gpsi_Q})$, les z{\'e}ros de~$\kappa_P$ sont aussi z{\'e}ros
de~$\frac{\kappa_P}{\kappa_P \circ \agot_P} - 1$. En bref,
 nous venons de montrer que si
$\frac{\kappa_P}{\kappa_P \circ \agot_P} - 1$
est non-nulle, elle a le m{\^e}me diviseur
que $\kappa_P$~: il
existe donc une constante~$c \in \CC$ telle que~:
$$
\frac{\kappa_P}{\kappa_P \circ \agot_P} - 1 = c \kappa_P
\qquad \text{ou encore~:} \qquad
\frac{1}{\kappa_P \circ \agot_P} = \frac{1}{\kappa_P} + c.
$$
Tenant compte de fait que~$\agot_P$ est d'ordre fini~$e$, on en d{\'e}duit que~$ce = 0$, puis
que~$c=0$, puis que~$\kappa_P \circ \agot_P = \kappa_P$. Autrement dit,~$\agot_P$ est un
automorphsime du rev{\^e}tement~$\kappa_P = g \circ \psi \circ \chi_P : \cF_P \to \PP^1$. Par
hypoth{\`e}se sur la fonction~$g$, il en r{\'e}sulte que~$\agot_P$ est la restriction {\`a}~$\cF_P$ d'un
{\'e}l{\'e}ment de~$A \times B$.

Quitte {\`a} composer~$\agot$ par l'inverse de cet {\'e}l{\'e}ment, on peut supposer que~$\agot_P$ est
l'identit{\'e}. En particulier,~$\agot$ fixe tous les points singuliers de~$\cK_0$ appartenant
{\`a}~$\cF_P$. Du coup,~$\agot$ stabilise n{\'e}cessairement toutes les composantes~$\cE_Q$ pour~$Q$
z{\'e}ro de~$g \circ \psi$. Par restriction, il induit donc un automorphisme~$\agot_Q$ de~$\cE_Q$
qui, d'apr{\`e}s~(\ref{eq:chiq}), est en fait un automorphisme de~$\chi_Q$. Comme de plus, il fixe un point
dans (et m{\^e}me toute) la fibre non ramifi{\'e}e au dessus de~$P$ du rev{\^e}tement
galoisien~$\chi_Q : \cE_Q \to \cB$, n{\'e}cessairement~$\agot_Q$ vaut l'identit{\'e}. Ainsi~$\agot$ fixe
point par point toutes les composantes~$\cE_Q$.

Il ne reste plus qu'{\`a} v{\'e}rifier qu'il en est de m{\^e}me des composantes~$\cF_{P'}$ pour~$P'$ un
autre z{\'e}ro de~$f$, i.e. distinct de~$P$. Notons~$\agot_{P'}$ la
restriction de~$\agot$ {\`a}~$\cF_{P'}$. 
On sait qu'elle co{\"\i}ncide avec la restriction d'un {\'e}l{\'e}ment de~$A \times B$ et qu'elle fixe
tous les points singuliers de~$\cK_0$ appartenant {\`a}~$\cF_{P'}$, c'est-{\`a}-dire tous les points
des fibres au dessus des z{\'e}ros de~$g \circ \psi$. Il suffit donc de v{\'e}rifier
que l'action de~$A \times B$ sur l'union de  ces fibres est libre. C'est assur{\'e}ment le cas pour les
{\'e}l{\'e}ments de~$B$ car les z{\'e}ros de~$g \circ \psi$ sont, par hypoth{\`e}se, non ramifi{\'e}s dans le
rev{\^e}tement galoisien~$\chi_{P'} : \cF_{P'} \to \cD$. \c Ca l'est encore pour les {\'e}l{\'e}ments
de~$A \times B$ car l'action de~$A$ sur les z{\'e}ros de~$g \circ \psi$ est aussi libre.
\end{preuve}

\subsection{D{\'e}formations}\label{subsection:deformation}

Nous allons maintenant d{\'e}former les deux courbes stables
pr{\'e}c{\'e}dentes pour faire en sorte qu'elles soient chacune la fibre
sp{\'e}ciale d'une famille de courbes.
Si $x\in K^s$ est un scalaire, il est naturel de consid{\'e}rer la courbe
$\cJ_x$ trac{\'e}e sur $\cI=(\cB\otimes_KL)\times \cD$ et d{\'e}finie par 
l'{\'e}quation $f(P)\times g(\psi(Q))=x$. On note $\cK_x$ l'image
r{\'e}ciproque de $\cJ_x$ par $\chi$.
Nous montrons dans ce paragraphe et dans le suivant que pour presque
tout $x$ dans $K$, la courbe $\cK_x$ est lisse et
 g{\'e}om{\'e}triquement  irr{\'e}ductible, 
de groupe d'automorphismes $A\times B$, et qu'elle a m{\^e}me corps des modules
et m{\^e}mes corps de d{\'e}finition que le rev{\^e}tement $\varphi$ de d{\'e}part.
Pour ce faire, nous voulons regarder  les familles
$(\cJ_x)_x$ et  $(\cK_x)_x$ comme
 des fibrations au-dessus  de $\PP^1$.
Mais
il faut {\^e}tre prudent : la famille de courbes $(\cJ_x)_x$ admet des points
bases. 
C'est pourquoi nous commen{\c c}ons par {\'e}clater la surface
$\cI=(\cB\otimes_KL)\times \cD$  le long de

$$c = ((f)_\infty\times (g\circ \psi)_0) \cup ((f)_0\times
(g\circ\psi)_\infty).$$

 Notons que $c$ est union  de 
$2\times \deg(f)\times \deg(g\circ\psi)$ points g{\'e}om{\'e}triques simples.
Soit $\cI_{\infty, \infty}\subset  \cI= (\cB\otimes_KL)\times \cD$ 
l'ouvert compl{\'e}mentaire de $((f)_\infty\times \cD) \cup (\cB\times
(g\circ\psi)_\infty)$ dans $(\cB\otimes_KL)\times \cD$.
On d{\'e}finit de m{\^e}me $\cI_{0,0}$, $\cI_{0,\infty}$, $\cI_{\infty,0}$.
Ces quatre ouverts recouvrent $(\cB\otimes_KL)\times \cD$.

On note $\PP^1_L=\Proj(L[X_0,X_1])$ la droite  projective sur $L$.
On note $F=1/f$ et $G=1/g$. 
Soit $\cJ_{\infty,0}\subset \cI_{\infty,0}\times 
\PP^1_L$ 
l'ensemble des $(P,Q,[X_0:X_1])$ tels que $f(P)X_0=G(\psi(Q))X_1$.
Soit $\cJ_{0,\infty}\subset \cI_{0,\infty}\times 
\PP^1_L$ 
l'ensemble des $(P,Q,[X_0:X_1])$ tels que $g(\psi(Q))X_0=F(P)X_1$.
Soit $\cJ_{\infty,\infty}\subset \cI_{\infty,\infty}\times 
\PP^1_L$ 
l'ensemble des $(P,Q,[X_0:X_1])$ tels que $f(P)g(\psi(Q))X_0=X_1$.
Soit $\cJ_{0,0}\subset \cI_{0,0}\times 
\PP^1_L$ 
l'ensemble des $(P,Q,[X_0:X_1])$ tels que $X_0=F(P)G(\psi(Q))X_1$.
On note $\cJ\subset \cI\times \PP^1_L$ le recollement de ces quatre ensembles,
$b_\cI : \cJ\rightarrow \cI$ la projection sur le premier facteur
et $j : \cJ\rightarrow \PP^1_L$ la projection
sur $\PP^1_L$. C'est un morphisme plat, projectif, surjectif.

Soit $\cK\subset \cS\times \PP^1_L$  l'image inverse de 
$\cJ$ par $\chi\times \Id$ o{\`u} $\Id : \PP^1_L\rightarrow \PP^1_L$ est l'identit{\'e}.
C'est l'{\'e}clatement de $\cJ$ le long de $\chi^{-1}(c)$. 
Notons que $\chi^{-1}(c)$ est union  de 
$\deg(\chi)\times \deg(f)\times \deg(g\circ\psi)$ points g{\'e}om{\'e}triques simples  car 
$\chi$ est non ramifi{\'e} au dessus de $c$.
En fait, $\cK$ est la normalisation de 
$\cJ$ dans $L(\cS)$. On note encore $\chi : \cK \rightarrow
\cJ$ le morphisme correspondant. 
Soit $b_\cS : \cK\rightarrow \cS$  la projection sur le 
premier facteur.
Soit $k : \cK \rightarrow \PP^1_L$
la projection sur le second facteur. C'est 
le morphisme compos{\'e} $k=j\circ \chi$. C'est un morphisme plat,
propre  et surjectif.


Soit  $\AA^1_L\subset \PP^1_L$ le spectre de $L[X]$ avec $X=X_1/X_0$.
La fonction $X$
permet d'identifier 
$\AA^1_L(K^s)$ avec $K^s$.
Pour tout point $x$ de $\AA^1(K^s)$  on note $\cJ_x$ 
la fibre de $j$ au-dessus de $x$ et $\cK_x$
la fibre de $k$ au dessus de $x$.  La restriction
de $b_\cI$ {\`a} $\cJ_x$ est une immersion ferm{\'e}e.
On peut donc voir $\cJ_x$ comme une courbe trac{\'e}e sur $\cI=
\cB\times\cD$.
De m{\^e}me, la restriction
de $b_\cS$ {\`a} $\cK_x$ est une immersion ferm{\'e}e.
On peut donc voir $\cK_x$ comme une courbe trac{\'e}e sur $\cS$.
On note $\cJ_\eta$ la fibre g{\'e}n{\'e}rique de $j$
et $\cK_\eta$ celle de $k$.
La fibre de $j$ en $0$ est  isomorphe
(par le morphisme $b_\cI$) {\`a} la courbe stable $\cJ_0$
introduite au paragraphe \ref{subsection:deux}.
De m{\^e}me, la fibre de $k$ en $0$ est  isomorphe
(par le morphisme $b_\cS$) {\`a} la courbe stable $\cK_0$
introduite au paragraphe \ref{subsection:deux}.

Montrons que  $\cJ_\eta/L(X)$
est  g{\'e}om{\'e}triquement connexe et que pour presque tout 
$x\in \AA^1(K^s)$ la courbe $\cJ_x$ est
g{\'e}om{\'e}tri\-que\-ment connexe  sur $L(x)$.
D'apr{\`e}s le th{\'e}or{\`e}me de factorisation de Stein \cite[Chapter 5,
Exercise 3.11]{Liu}, on peut factoriser $j : \cJ/L
\rightarrow \AA^1_L$ en $j_f\circ j_c$ avec $j_c$  
{\`a} fibres g{\'e}om{\'e}triquement
connexes et $j_f$ fini et dominant.
La fibre de $j_f$ au dessus de $0$ est triviale car $\cJ_0$ est connexe 
et r{\'e}duite. Donc le degr{\'e} de $j_f$ est $1$ 
d'apr{\`e}s \cite[Chapter
5,Exercise 1.25]{Liu}. Donc $j_f$ est un isomorphisme  au
dessus d'un ouvert  
de $\AA^1_L$. La fibre g{\'e}n{\'e}rique $\cJ_\eta$ est g{\'e}om{\'e}triquement
connexe sur $L(X)$ et 
pour presque tout $x\in \AA^1(K^s)$ la courbe $\cJ_x$ est
g{\'e}om{\'e}tri\-que\-ment connexe  sur $L(x)$.
Le m{\^e}me raisonnement montre que 
$\cK_\eta$ est g{\'e}om{\'e}triquement
connexe sur $L(X)$ et 
pour presque tout $x\in \AA^1(K^s)$ 
la courbe $\cK_x$ est g{\'e}om{\'e}tri\-que\-ment connexe sur $L(x)$.

Montrons maintenant
que $\cJ_\eta$ est lisse. 
Elle  est lisse en dehors des
points~$(P,Q) \in \cJ_\eta\subset \cB \times \cD$
o{\`u}~$df(P) =0$ et $d(g \circ \psi)(Q) = 0$. 
De tels points sont d{\'e}finis sur $K^s$ donc la fonction
$f(P)\times g(\psi(Q))$ ne peut pas prendre la valeur 
transcendante 
$X$ en ces points.

Le lieu $\Gamma$ de ramification de $\chi$ coupe
transversalement la fibre $\cJ_0$. Donc il coupe
transversalement la  fibre g{\'e}n{\'e}rique $\cJ_\eta$.
Donc  $\cK_\eta$ est lisse elle aussi.
Ainsi, pour presque tout $x\in K^s$ les fibres $\cJ_x$ et $\cK_x$
sont lisses.

Enfin, la connaissance de~$\Aut_{K^s}^{\text{d{\'e}f.}} (\mathcal{K}_0)$
permet
de montrer que~$\Aut_{\overline{K(X)}}(\mathcal{K}_\eta) \simeq A \times B$.
Soit $R=K^s[[X]]$.
La courbe  localis{\'e}e $\hat \cK = \cK \otimes_{\AA^1_L}\Spec (R)$ est 
stable sur 
le spectre de $R$. D'apr{\`e}s \cite[Chapter 10,
Proposition 3.38, Remarque 3.39]{Liu} le foncteur 
automorphismes $t\mapsto \Aut _t(\hat \cK_t)$ est repr{\'e}sentable 
par un sch{\'e}ma fini
et non-ramifi{\'e} sur $\Spec  R$ et l'application de
sp{\'e}cialisation   $\Aut_{K^s((X))}(\hat \cK_\eta)
\rightarrow \Aut_{K^s} (\cK_0)$ est injective. Elle prend ses valeurs dans le
sous-groupe
des $K^s$-automorphismes d{\'e}formables de $\cK_0$ d'apr{\`e}s 
le lemme \ref{lem:defo}. 
Comme $\Spec  R$ n'a pas de rev{\^e}tement
non-ramifi{\'e}, on en d{\'e}duit que le groupe des automorphismes
g{\'e}om{\'e}triques de la fibre g{\'e}n{\'e}rique satisfait
$$
A\times B\subset \Aut_{\overline{K(X)}}(\cK_\eta) \subset 
\Aut_{{K^s((X))}}(\cK_\eta)\subset \Aut^{\text{d{\'e}f.}}_{K^s}(\cK_0).
$$
Comme le dernier groupe, on l'a vu, est {\'e}gal {\`a}~$A \times B$, on en d{\'e}duit bien
que~$\Aut_{\overline{K(X)}}(\cK_\eta) = A \times B$.

\subsection{Corps de modules, corps de d{\'e}finition dans les fibres de~$\mathcal{K}$}

Pour conclure nous montrons l'existence de 
fibres~$\mathcal{K}_x$ avec~$x \in \AA^1(K)$, lisses, g{\'e}om{\'e}trique\-ment
irr{\'e}ductibles, 
ayant m{\^e}me corps des modules que~$\mathcal{S}$ (c'est-{\`a}-dire~$K$ par hypoth{\`e}se),
et m{\^e}mes corps de d{\'e}finition.

On vient de voir
 que pour presque tout~$x \in \AA^1(K)$, la fibre~$\mathcal{K}_x$ est 
lisse et g{\'e}om{\'e}triquement
connexe, donc g{\'e}om{\'e}triquement irr{\'e}ductible. Compte tenu  du lemme
\ref{lem:specialisation} sur la sp{\'e}cialisation du groupe
des automorphismes dans une famille de courbes rationnelles,
pour presque tout~$x \in \AA^1(K)$, le
groupe de  $K^s$-automorphismes  de la fibre
$\cK_x$ est isomorphe au  groupe 
des $\overline{K(X)}$-automorphismes
de  la fibre
g{\'e}n{\'e}rique.
Comme le groupe d'automorphismes 
 de la fibre g{\'e}n{\'e}rique est aussi isomorphe
au groupe $A\times B$ 
des automorphismes de~$\mathcal{S}$, on en d{\'e}duit que la restriction 
induit, pour presque tout $x$,  un
isomorphisme~:
$$
\Aut_{K^s}(\mathcal{S}) \overset{\simeq}{\longrightarrow} \Aut_{K^s}(\mathcal{K}_x).
$$

Soit donc $x\in K$ tel que $\cK_x$ soit lisse et  g{\'e}om{\'e}triquement  irr{\'e}ductible et
$\Aut_{K^s}(\cK_x)=A\times B$.
Montrons que $\mathcal{K}_x$ poss{\`e}de toutes les propri{\'e}t{\'e}s de module et de d{\'e}finition
escompt{\'e}es.

Tout d'abord, la construction de~$\mathcal{K}$ n'impliquant pas d'extension
des scalaires, {\`a} tout mod{\`e}le de $\varphi$ d{\'e}fini sur un corps
$M\supset K$ on associe un $M$-mod{\`e}le  de la courbe~$\mathcal{K}_x$.

Ensuite, la surface~$\mathcal{S}$ a~$K$ pour corps de modules, et pour
tout~$\sigma \in \Gal(K)$, on sait qu'il existe un
$K^s$-isomorphisme~$\varphi_\sigma : \mathcal{S} \to {}^\sigma \mathcal{S}$. Par restriction,
ces~$\varphi_\sigma$ induisent des $K^s$-isomorphismes entre~$\mathcal{K}_x$
et~$^{\sigma}\mathcal{K}_x$.
Donc le  corps des modules de $\cK_x$  est $K$.

En fait, tout $K^s$-isomorphisme entre $\cK_x$ et $^{\sigma}\mathcal{K}_x$
est restriction d'un $K^s$-isomorphisme entre $\cS$ et $^\sigma\cS$. Il
existe en effet 
un tel isomorphisme, et les groupes d'automorphismes de $\cK_x$ et de $\cS$
sont les m{\^e}mes. Comme la correspondance est bijective, 
toute donn{\'e}e de
descente pour $\cK_x$ s'{\'e}tend 
en une donn{\'e}e de descente pour $\cS$.

Cela termine la d{\'e}monstration du th{\'e}or{\`e}me
\ref{theoreme:courbe}.
Le th{\'e}or{\`e}me \ref{theoreme:courbe} permet de d{\'e}montrer
le th{\'e}or{\`e}me \ref{theoreme:courbex} {\`a} l'aide du th{\'e}or{\`e}me
\ref{theoreme:catrevetement}.
Notons que Jean-Fran{\c c}ois Mestre a construit dans~\cite{Mestre} des obstruction
{\it locales} {\`a} l'existence d'un mod{\`e}le d{\'e}fini sur le corps des modules d'une courbe.

\section{Quelques r{\'e}sultats interm{\'e}diaires d'ordre g{\'e}n{\'e}ral}

\subsection{Sur les courbes et les produits de deux courbes}

\begin{lem}\label{lem:adjonction}
Soit $K$ un corps de caract{\'e}ristique nulle.
Soient  $\cB$ et $\cC$  deux $K$-courbes alg{\'e}briques projectives,
r{\'e}duites, lisses et g{\'e}om{\'e}triquement
irr{\'e}\-duc\-tibles  de genres $\beta$ et $\gamma$ respectivement.
On fixe $P$ un point g{\'e}om{\'e}trique
  de $\cB$ et $Q$ un point
g{\'e}om{\'e}trique  de $\cC$ et on note abusivement $\cB=\cB\times Q$ et $\cC=P\times \cC$.
Soit  $\Gamma$ un diviseur de $\cB\times \cC$ de 
bidegr{\'e} $(b,c)$, i.e. tel que~$b=\cB\cdot\Gamma$ et $c=\cC\cdot\Gamma$. Le genre arithm{\'e}tique
virtuel $\pi$ de $\Gamma$ est au plus $1+bc+c(\beta-1)+b(\gamma-1)$. Dans le cas $b=c$ cette borne
vaut $1+2b(\beta -1)+b^2$.
\end{lem}

\begin{preuve}
On s'inspire ici de la preuve de Weil de l'hypoth{\`e}se de Riemann pour les courbes
(cf.~\cite[Exercice~V-1.10]{Hartshorne}).

La classe d'{\'e}quivalence alg{\'e}brique du diviseur canonique de $\cB\times\cC$
est $K=2(\beta -1)\cC+2(\gamma -1)\cB$. Gr{\^a}ce {\`a} la formule
d'adjonction, on sait que~$\pi=\frac{D\cdot(D+K)}{2}+1$ (cf.~\cite[Exercice~V-1.3]{Hartshorne}).
On a donc~:
$$
\pi = \frac{D\cdot\left(D + 2(\beta -1)\cC+2(\gamma -1)\cB\right)}{2} + 1
    = \frac{D \cdot D + 2c(\beta - 1) + 2b(\gamma - 1)}{2} + 1,
$$
et il ne reste plus qu'{\`a} majorer l'auto-intersection~$D \cdot D$. L'in{\'e}galit{\'e} de
Castelnuovo et Severi (cf.\cite[Exercice~V-1.9]{Hartshorne}) nous apprend
que~$D \cdot D \leq 2bc$, d'o{\`u} le r{\'e}sultat.
\end{preuve}

\begin{lem} \label{lem:fibrescoupenttransversalement}
Soit $K$ un corps de caract{\'e}ristique nulle.
Soient~$\mathcal{B}$ et~$\mathcal{C}$ deux $K$-courbes 
algebriques projectives lisses,  g{\'e}om{\'e}triquement irr{\'e}ductibles et
r{\'e}duites.
Soit $\Gamma$ un diviseur
sans multiplicit{\'e} de la surface~$\mathcal{B} \times \mathcal{C}$. 
Soit $f \in K(\mathcal{B})$ une fonction
non constante. Pour  tout 
{\'e}l{\'e}ment $x$ de $K^s$ sauf un nombre
fini,  le
diviseur~$(f)_x \times \mathcal{C}$ coupe transversalement~$\Gamma$, 
o{\`u}~$(f)_x$ d{\'e}signe le diviseur
$f^{-1}(x)$.
\end{lem}

\begin{preuve}
Notons~$p_\cB : \cB \times \cC \to \cB$ la premi{\`e}re projection. L'ensemble~$E$ des points
de~$\cB(K^s)$ tels que dans~$p_\cB^{-1}(P)$ il y ait soit un point singulier
de~$\Gamma$, soit un point ramifi{\'e} du morphisme~$\Gamma \to \cB$ obtenu
par restriction de~$p_\cB$, 
ou tels que la fibre~$p_\cB^{-1}(P)$ soit contenue dans~$\Gamma$,
 est fini. Donc  pour  tout $x\in K^s$ sauf un nombre fini,
 la fibre $f^{-1}(x)$
{\'e}vite~$E$ et elle est simple. 
\end{preuve}

\begin{lem}\label{lem:autofibres}
Soit $K$ un corps de caract{\'e}ristique nulle.
Soit $\cB$ une  $K$-courbe alg{\'e}brique projective
lisse, g{\'e}om{\'e}triquement irr{\'e}ductible et r{\'e}duite. On suppose que
le genre de $\cB$ est au moins $2$.
Soit $f\in K(\mathcal{B})$ une fonction  non-constante.
On note $G$ le groupe des $K^s$-automorphismes de $f$.
C'est l'ensemble des $K^s$-automorphismes $\agot$ 
de $\cB$  tels que $f \circ \agot = f$. Pour
tout $x\in \PU (K^s)$,  on note 
$(f)_x=f^{-1}(x)$ la fibre au dessus de $x$
et $G_x$ le groupe des $K^s$-automorphismes
de $\cB$ qui stabilisent 
l'ensemble  des $K^s$-points de $(f)_x$.

Pour tout $x$ de $\PU (K^s)$ sauf un nombre fini on a $G_x=G$.
\end{lem}

\begin{preuve}
Le  groupe $H=\Aut_{K^s}(\cB)$ des $K^s$-automorphismes
de $\cB$ est fini car le genre de $\cB$ est au moins {\'e}gal {\`a} $2$.
Soit $\agot$ un automorphisme dans $H-G$ et 
soit  $x\in \PU (K^s)$. Supposons
que  les $K^s$-points de $(f)_x$ sont stabilis{\'e}s
par $\agot$. Soit $P$ l'un d'eux. Alors $f\circ \agot (P)=f(P)=x$.
Ainsi $P$ est un z{\'e}ro de la fonction non-nulle $f\circ \agot -f$.
Pour chaque $\agot$ on trouve un nombre fini de tels z{\'e}ros.
Et les $\agot$ sont en nombre fini. Donc les $x$, qui sont
images par $f$ de tels $P$, sont eux aussi en nombre
fini.
\end{preuve}

\begin{lem}\label{lem:fonctionsauto}
Soit $K$ un corps de caract{\'e}ristique nulle. 
Soit $\cB$ une    $K$-courbe alg{\'e}brique projective
lisse, g{\'e}om{\'e}triquement irr{\'e}ductible et r{\'e}duite. 
Soit $L\subset K^s$ une extension alg{\'e}brique de $K$.
Soit $\cC$ une    $L$-courbe alg{\'e}brique projective
lisse, g{\'e}om{\'e}triquement irr{\'e}ductible et r{\'e}duite. On suppose que
le genre de $\cB$ est au moins $2$.
Soit $\varphi : \cC\rightarrow \cB\otimes_K L$ un $L$-morphisme non-constant.
On note $G$ le groupe des $K^s$-automorphismes de $\varphi$.
C'est l'ensemble des $K^s$-automorphismes $\agot$ 
de $\cC$  tels que $\varphi \circ \agot = \varphi$. 
Pour toute fonction $f$ dans $K^s(\cB)$ on note 
$G_f=\Aut_{K^s}(f\circ \varphi)$
le groupe des $K^s$-automorphismes de $f\circ \varphi$.
On a $G\subset G_f$.
On note $V\subset K(\cB)$ l'ensemble des fonctions $f$ 
telles que $G_f\not = G$.
L'ensemble $V$ est contenu dans une union finie de $K$-alg{\`e}bres
 strictes de $K(\cB)$.
\end{lem}


\begin{preuve}
Par l'{\'e}quivalence cat{\'e}gorielle entre les corps de fonctions et les courbes
(projectives lisses, etc \ldots), l'{\'e}nonc{\'e} {\`a} prouver se ram{\`e}ne {\`a} un r{\'e}sultat concernant les corps de
fonctions~$K^s(f) \subset K^s(\cB) \subset K^s(\cC)$ et les groupes qui 
entrent en jeu sont les
suivants~:
$$
\begin{cases}
G = \Aut_{K^s}(\varphi) = \Aut_{K^s(\cB)}(K^s(\cC)), \\
G_f = \Aut_{K^s}(f \circ \varphi) = \Aut_{K^s(f)}(K^s(\cC)), \\
A = \Aut_{K^s}(\cC) = \Aut_{K^s}(K^s(\cC)),
\end{cases}
\quad\Rightarrow\quad
G \subset G_f \subset A.
$$
Dans ce contexte l'ensemble~$V$ admet la description~:
$$
V
=
\left(
\bigcup_{\agot \in A \setminus G}
K^s(\cC)^\agot \cap K^s(\cB)
\right) \cap K(\cB)
=
\bigcup_{\agot \in A \setminus G}
K^s(\cC)^\agot \cap K(\cB).
$$
Cette r{\'e}union est finie, le groupe~$A$ l'{\'e}tant puisque~$\cC$ est de genre~$\geq 2$.
De plus, comme~$\agot \not \in G$, chacune des sous-extensions~$K^s(\cC)^\agot \cap K^s(\cB)$ est
une sous-extension stricte de~$K^s(\cB)$, contenant~$K^s$~; en d{e}scendant {\`a}~$K$, on pr{\'e}serve
l'inclusion stricte~$K^s(\cC)^\agot \cap K(\cB) \subsetneq K(\cB)$.
\end{preuve}
\subsection{D{\'e}formation d'automorphismes}

Dans ce paragraphe nous donnons une condition n{\'e}cessaire
pour qu'un automorphisme d'une courbe nodale s'{\'e}tende {\`a} une d{\'e}formation
de cette courbe.

Soit~$R$ un anneau de valuation discr{\`e}te complet,~$\pi$ une uniformisante et~$k$ son corps r{\'e}siduel
suppos{\'e} alg{\'e}briquement clos. On consid{\`e}re~$\cK$ une courbe stable sur~$\Spec(R)$~; la fibre
g{\'e}n{\'e}rique est not{\'e}e~$\cK_\eta$, la fibre sp{\'e}ciale~$\cK_0$. On consid{\`e}re~$T$ un point singulier
de~$\cK_0$. D'apr{\`e}s~\cite[Chap~10, Corollaire~3.22]{Liu}, on sait que l'anneau local compl{\'e}t{\'e}
de~$\cK$ en~$T$ est de la forme~:
$$
\hat{\cO}_{\cK, T} = R[[f,g]]/\langle f g - \pi^e \rangle
$$
avec~$e \in \NN$. Cet entier s'appelle {\em l'{\'e}paisseur} de~$\cK$ en~$T$. On dit aussi que~$f$
et~$g$ forment un syst{\`e}me de coordonn{\'e}es de~$\cK$ en~$T$. En r{\'e}duisant
modulo~$\pi$, on obtient l'anneau local compl{\'e}t{\'e} de~$\cK_0$ en~$T$~:
$$
\hat{\cO}_{\cK_0, T} = \hat{\cO}_{\cK, T} / \langle \pi \rangle
= k[[\overline{f},\overline{g}]]/\langle \overline{f} \overline{g} \rangle,
$$
o{\`u}~$\overline{f} = f \bmod{\pi}$ et~$\overline{g} = g \bmod{\pi}$.
Comme en tout point double ordinaire, on sait que, localement en~$T$,~$\cK_0$ poss{\`e}de deux
branches~$\cF$ et~$\cG$ (composantes irr{\'e}ductibles)~; les fonctions~$\overline{f}$
et~$\overline{g}$ en sont les uniformisantes respectives. On note~$P$ et~$Q$ les points de~$\cF$
et~$\cG$ correspondant {\`a}~$T$.

Soit~$T'$ un autre point point singulier de~$\cK_0$ auxquel on associe les
donn{\'e}es~$f',g',e',\cF', \cG'$ comme pr{\'e}c{\'e}demment.

Consid{\'e}rons~$\agot$ un automorphime de~$\cK$ sur~$R$ tel que~$\agot(T) = T'$
et~$\agot(\cF) = \cF'$, ~$\agot(\cG) = \cG'$. On v{\'e}rifie ais{\'e}ment que les
fonctions~$f' \circ \agot$ et~$g' \circ \agot$ forment aussi un syst{\`e}me de coordonn{\'e}es de~$\cK$
en~$T$. On en d{\'e}duit que~$e' = e$ mais aussi que~$f' \circ \agot / f$ et~$g' \circ \agot / g$
sont des inversibles de~$\hat{\cO}_{\cK, T}$ ({\`a} chaque fois le num{\'e}rateur et le d{\'e}nominateur ont
m{\^e}me diviseur de Weil). Comme de plus on
a~$f \times g = \pi^e = f' \circ \agot \times g' \circ \agot$, il en r{\'e}sulte
que~$\frac{f' \circ \agot}{f}(T) \times \frac{g' \circ \agot}{g}(T) = 1$.
En r{\'e}duisant modulo~$\pi$ et en se pla{\c c}ant sur chacune des deux branches locales de~$\cK_0$
en~$T$, on en d{\'e}duit que~:

\begin{equation}\label{eq:deformablegene}
\frac{\overline{f}' \circ \overline{\agot}}{\overline{f}}(P) \times \frac{\overline{g}' \circ \overline{\agot}}{\overline{g}}(Q) = 1.
\end{equation}

\begin{definition}\label{defi:deform}
Soit $R$ un anneau de valuation discr{\`e}te complet
de corps r{\'e}siduel $k$
alg{\'e}brique\-ment clos. Soit $\cK$ une courbe
semi-stable sur $\Spec(R)$. La fibre g{\'e}n{\'e}rique 
de $\cK$ est suppos{\'e}e lisse. On se donne un syst{\`e}me
de coordonn{\'e}es en chaque point singulier de la fibre sp{\'e}ciale $\cK_0$.
Soit $\bagot$ un automorphisme de 
$\cK_0/k$. On dit que $\bagot$ est {\it d{\'e}formable} dans $\cK/\Spec(R)$
si
pour tout point singulier $T$ de $\cK_0$, l'image $\bagot(T)$
a la m{\^e}me {\'e}paisseur que $T$ dans $\cK$,
  et si l'{\'e}galit{\'e} (\ref{eq:deformablegene})
est satisfaite.
\end{definition}

\begin{lem}\label{lem:defo}
Avec les notations de la d{\'e}finition \ref{defi:deform}, l'ensemble
des automorphismes de $\cK_0/k$  d{\'e}\-for\-mables dans $\cK/\Spec(R)$  
est un sous-groupe de $\Aut_k(\cK_0)$.
Si $\agot$ est un automorphisme de $\cK$ sur $\Spec(R)$, sa r{\'e}duction
$\bagot = \agot \bmod \pi$ est un automorphisme de $\cK_0/k$
d{\'e}formable dans $\cK/\Spec(R)$.
\end{lem}

On pourra comparer cet {\'e}nonc\'e {\`a} \cite[Theorem 3.1.1]{Wewers}
qui concerne les d{\'e}formations de morphismes entre deux courbes
distinctes.

\subsection{Sur les automorphismes d'une famille de courbes}
Nous d{\'e}montrons dans cette section un  lemme de sp{\'e}cialisation du groupe
des automorphismes d'une courbe d{\'e}pendant d'un param{\`e}tre.

Soit $K$ un corps de caract{\'e}ristique nulle  et 
soit $\cU$  une courbe affine  lisse,
g{\'e}om{\'e}triquement irreductible et r{\'e}duite sur
$K$. 
Soit  $X$ une surface quasi-projective lisse, g{\'e}om{\'e}triquement irr{\'e}duct\-ible et r{\'e}duite. Soit 
$\pi : X \rightarrow \cU $  un morphisme  projectif, lisse, surjectif,  de dimension
relative $1$. On suppose que  pour tout point 
g{\'e}om{\'e}trique $\epsilon\in \cU(K^s)$
de $\cU$,  la fibre $X_\epsilon$
 est  g{\'e}om{\'e}triquement irr{\'e}ductible.

Toutes les fibres ont le m{\^e}me genre g{\'e}om{\'e}trique $g$ suppos{\'e} au
moins {\'e}gal {\`a} $3$.

On veut montrer qu'il existe un $K$-ouvert non-vide $\cV$ de $\cU$
tel que pour tout point
g{\'e}om{\'e}trique  $\epsilon \in \cV(K^s)$ le groupe des automorphismes sur $K^s$
de la fibre en $\epsilon$ soit {\'e}gal au groupe 

$$\Aut_{\overline{K(\cU)}}(X_\eta\otimes_{K(\cU)}\overline{K(\cU)})$$
\noindent  des automorphismes 
de la fibre g{\'e}n{\'e}rique $X_\eta$.

Supposons
d'abord que  la fibre g{\'e}n{\'e}rique $X_\eta$ 
est non-hyperelliptique (ce qui revient {\`a} dire
que le nombre de points de Weierstrass est plus grand que $2g+2$).

Quitte {\`a} remplacer
$K$ par une extension finie et $\cU$ par
un rev{\^e}tement fini {\'e}tale d'un ouvert non-vide de 
$\cU$, on  peut supposer
que les points de Weierstrass de la fibre g{\'e}n{\'e}rique
sont d{\'e}finis sur le corps $K(\cU)$. L'adh{\'e}rence de Zariski de 
ces points
d{\'e}finit des sections horizontales de la famille
$\pi : X\rightarrow \cU$. 
Quitte {\`a} restreindre $\cU$  un peu plus, on peut supposer que 
ces sections ne se croisent pas. Donc aucune fibre de la famille
n'est hyperelliptique.

Il existe $g$ points de Weierstrass $P_{1,\eta}$, $P_{2,\eta}$, \ldots , 
$P_{g,\eta}$ de la fibre
g{\'e}n{\'e}rique qui sont lin{\'e}airement ind{\'e}pendants dans son
mod{\`e}le canonique. 
On choisit une base $(\omega_1,
\ldots , \omega_g)$
des diff{\'e}rentielles holomorphes de la fibre g{\'e}n{\'e}rique, adapt{\'e}e {\`a} ces $g$
points (i.e. $\omega_i(P_{j,\eta})=0$ si $i\not = j$). 
Quitte {\`a} restreindre un
peu plus l'ouvert $\cU$ on peut supposer que ces $g$ formes se sp{\'e}cialisent
en tout point $x\in \cU$ pour former une base  des formes holomorphes de la
fibre
en $x$.
  
On obtient donc une famille de plongements canoniques

\begin{equation*}
\xymatrix{
X   \,\,\, \ar@{^{(}->}[r] \ar@{->}[d] &\cU\times \PP^{g-1} \ar@{->}[dl] \\
\cU&\\  }  \label{diagramme:plong}
\end{equation*}

Pour chaque fibre de la   famille  $\pi : X \rightarrow \cU$, le groupe d'automorphismes peut-{\^e}tre vu comme
sous-groupe du groupe des permutations $S_W$ des points de Weierstrass\footnote{Donc tous les automorphismes de $X_\eta$ sont d{\'e}finis sur $K^s(\cU)$,
c'est-{\`a}-dire $\Aut_{\overline{K(\cU)}}(X_\eta)=
\Aut_{K^s(\cU)}(X_\eta)$.}.
En particulier, si $\epsilon \in \cU(K^s)$ est un point ferm{\'e} et $\eta $ le point g{\'e}n{\'e}rique
alors $\Aut_{K^s(\cU)}(X_\eta)\subset \Aut_{K^s}(X_\epsilon )\subset
S_W$.

On veut montrer que les $\epsilon$ tels que la premi{\`e}re
inclusion soit stricte sont
contenus
dans un ferm{\'e} strict de $\cU$. 
On choisit une permutation $\sigma \in S_W$ telle que 
$\sigma \not \in \Aut_{K^s(\cU)}(X_\eta)$. 
On associe {\`a}  $\sigma$ une transformation projective lin{\'e}aire de $\PP^{g-1}$ par
prolongement de l'action lin{\'e}aire sur les $P_1$, \ldots , $P_g$
(adh{\'e}rences de Zariski des $P_{1,\eta}$, $P_{2,\eta}$, \ldots,
$P_{g,\eta}$). Ce prolongement
est unique car les $P_i$ sont ind{\'e}pendants. On le note $\tau$. C'est un
automorphisme
de $\PP^{g-1}$ sur $\cU$.
Quitte {\`a} remplacer
$K$ par une extension finie et $\cU$ par
un rev{\^e}tement fini {\'e}tale d'un ouvert non-vide de 
$\cU$, on  peut supposer
qu'il existe un $K(\cU)$-point $Q_\eta\in X_\eta$ de 
la fibre g{\'e}n{\'e}rique tel que  $\tau (Q_\eta)$ ne soit
pas dans cette fibre. On prolonge ce point en une section rationnelle $Q$
de $\pi$.
L'ensemble des $x\in \cU$ tels que $\tau (Q_x)$ soit dans la fibre $X_x$ est
un ferm{\'e} strict de $\cU$.

Dans le cas o{\`u} la fibre g{\'e}n{\'e}rique est
hyperelliptique, notons  $c_\eta : X_\eta\rightarrow \PP^{1}/K(\cU)$
le morphisme canonique de degr{\'e} $2$, quotient par l'involution.
Soient $P_{1,\eta}$, \ldots, $P_{2g+2,\eta}$  les points de $\PP^{1}$ o{\`u}
$c_\eta$ se ramifie. Quitte {\`a} remplacer $K$ par une extension finie
et $\cU$ par un rev{\^e}tement {\'e}tale d'un ouvert non-vide,
on peut supposer 
que
le morphisme $c_\eta$ s'\'etend en un morphisme $c : X/\cU
\rightarrow  \PP^1/\cU$, que tous  les $P_{i,\eta}$
 sont d{\'e}finis sur $K(\cU)$, et 
que leurs  adh{\'e}rences
de Zariski  dans $\PP^1/\cU$ (not{\'e}es $P_i$)
ne se croisent pas.

Soit $x$ un point de $\cU$.
Tout $\overline{L(x)}$-automorphisme de $X_x$ normalise l'involution et induit un 
automorphisme de $\PP^1/\overline{L(x)}$ qui stabilise les $P_{i,x}$. 
On note $G_x$ l'ensemble de ces automorphismes.
Ce groupe est isomorphe
au  quotient du
groupe des $\overline{L(x)}$-automorphismes de $X_x$ par l'involution.

On note $\cM$ l'ensemble des injections de $\{1,2,3\}$ dans
$\{1,2,\ldots, 2g+2\}$.
Le groupe $G_x$  peut-{\^e}tre vu comme
sous-ensemble de  $\cM$ : on associe {\`a} $1$ (resp. $2$, $3$) l'indice de l'image
de $P_{1,x}$ (resp. $P_{2,x}$, $P_{3,x}$). 
En particulier
si $\epsilon \in \cU(K^s)$ est un point ferm{\'e} et $\eta $ le point g{\'e}n{\'e}rique
alors $G_\eta \subset G_\epsilon \subset \cM$.

On veut montrer que les $\epsilon$ tels que l'inclusion soit stricte sont
contenus
dans un ferm{\'e} strict de $\cU$. 
On choisit une injection  $\sigma \in \cM$ telle que 
$\sigma \not \in G_\eta$. 
On associe {\`a}  $\sigma$ la transformation projective 
lin{\'e}aire de $\PP^1/\cU$ qui envoie $P_i$ sur $P_{\sigma(i)}$
pour $i\in \{1,2,3\}$.  On la note $\tau$. C'est un
automorphisme
de $\PP^1$ sur $\cU$.
Il existe un indice $j$ entre $4$ et $2g+2$ tel que
pour tout $i\in \{1,2, \ldots, 2g+2\}$ 
$\tau(P_{j,\eta})\not =P_{i,\eta}$.
L'ensemble de $x\in \cU$ tels que $\tau (P_{j,x})$ soit 
{\'e}gal {\`a} l'un des $P_{i,x}$ est
un ferm{\'e} strict de $\cU$.

\begin{lem}\label{lem:specialisation}
Soit $K$ un corps de caract{\'e}ristique nulle  et $\cU$ une
courbe affine lisse g{\'e}om{\'e}trique\-ment  irr{\'e}ductible et r{\'e}duite
sur $K$. Soit
$X$ une surface quasi-projective lisse
g{\'e}om{\'e}triquement irr{\'e}ducti\-ble et  r{\'e}duite. Soit
$\pi : X \rightarrow \cU $  un
morphisme surjectif, projectif et lisse de dimension
relative $1$. On suppose que pour tout point  $\epsilon$
de $\cU$, la fibre $X_\epsilon$ en $\epsilon$
  est g{\'e}om{\'e}triquement irr{\'e}ductible.
On note $\eta$ le point g{\'e}n{\'e}rique de $\cU$ et
$\bar S_\eta = S_\eta\otimes_{K(\cU)}\overline{K(\cU)}$ 
la fibre g{\'e}n{\'e}rique vue comme courbe sur la cl{\^o}ture alg{\'e}brique du corps des
fonctions
de la base $\cU$.

Il existe un $K$-ouvert non-vide $\cV$ de $\cU$ tel que pour tout point
ferm{\'e} $\epsilon \in \cV(K^s)$ le groupe des automorphismes sur $K^s$
de la fibre en $\epsilon$ soit {\'e}gal au groupe 
$\Aut_{\overline{K(\cU)}}(\bar S_\eta )$
\noindent  des automorphismes 
de $\bar S_\eta$.
\end{lem}

\bibliographystyle{alpha}

\begin{thebibliography}{Wew99}

\bibitem[CR04]{roscouveignes}
Jean-Marc Couveignes and Nicolas Ros.
\newblock Des obstructions globales à la descente des revêtements.
\newblock {\em Acta Arithmetica}, 114(4):331--348, 2004.

\bibitem[DD99]{DD}
J.-C. Douai and P.~Dèbes.
\newblock Gerbes and covers.
\newblock {\em Comm. in Algebra}, 27/2:577--594, 1999.

\bibitem[Gir71]{giraud}
J.~Giraud.
\newblock {\em Cohomologie non-abélienne}, volume 179 of {\em Grundlehren Math.
  Wiss.}
\newblock Springer, 1971.

\bibitem[Har77]{Hartshorne}
Robin Hartshorne.
\newblock {\em Algebraic Geometry}, volume~52 of {\em Graduate Texts in
  Mathematics}.
\newblock Springer, 1977.

\bibitem[Liu02]{Liu}
Qing Liu.
\newblock {\em Algebraic Geometry and Arithmetic Curves}, volume~6 of {\em
  Oxford Graduate Texts in Mathematics}.
\newblock Oxford, 2002.

\bibitem[Mes91]{Mestre}
Jean-François Mestre.
\newblock Construction de courbes de genre~$2$ à partir de leurs modules.
\newblock In Teo Mora and Carlo Traverso, editors, {\em Effective Methods in
  Algebraic Geometry}, volume~94 of {\em Progress in Mathematics}, pages
  313--334. Birkäuser, 1991.

\bibitem[Wew99]{Wewers}
Stefan Wewers.
\newblock Deformation of tame admissible covers of curves.
\newblock In Helmut {\sc Völkein}, David {\sc Harbater}, Peter {\sc Müller},
  and J.G. {\sc Thompson}, editors, {\em Aspects of Galois theory}, volume 256
  of {\em London Mathematical Society Lecture Note Series}, pages 239--282.
  Cambridge, 1999.

\end{thebibliography}

\end{document}